\documentclass[12pt]{article}
\usepackage{latexsym, amssymb}
\DeclareMathAlphabet\gothic{U}{euf}{m}{n}

\makeatletter
\def\eqnarray{\stepcounter{equation}\let\@currentlabel=\theequation
\global\@eqnswtrue
\tabskip\@centering\let\\=\@eqncr
$$\halign to \displaywidth\bgroup\hfil\global\@eqcnt\z@
  $\displaystyle\tabskip\z@{##}$&\global\@eqcnt\@ne
  \hfil$\displaystyle{{}##{}}$\hfil
  &\global\@eqcnt\tw@ $\displaystyle{##}$\hfil
  \tabskip\@centering&\llap{##}\tabskip\z@\cr}

\def\endeqnarray{\@@eqncr\egroup
      \global\advance\c@equation\m@ne$$\global\@ignoretrue}

\def\@yeqncr{\@ifnextchar [{\@xeqncr}{\@xeqncr[5pt]}}
\makeatother

\begin{document}
\bibliographystyle{tom}

\newtheorem{lemma}{Lemma}[section]
\newtheorem{thm}[lemma]{Theorem}
\newtheorem{cor}[lemma]{Corollary}
\newtheorem{voorb}[lemma]{Example}
\newtheorem{rem}[lemma]{Remark}
\newtheorem{prop}[lemma]{Proposition}
\newtheorem{stat}[lemma]{{\hspace{-5pt}}}
\newtheorem{obs}[lemma]{Observation}
\newtheorem{defin}[lemma]{Definition}

\newenvironment{remarkn}{\begin{rem} \rm}{\end{rem}}
\newenvironment{exam}{\begin{voorb} \rm}{\end{voorb}}
\newenvironment{defn}{\begin{defin} \rm}{\end{defin}}
\newenvironment{obsn}{\begin{obs} \rm}{\end{obs}}

\newcommand{\gota}{\gothic{a}}
\newcommand{\gotb}{\gothic{b}}
\newcommand{\gotc}{\gothic{c}}
\newcommand{\gote}{\gothic{e}}
\newcommand{\gotf}{\gothic{f}}
\newcommand{\gotg}{\gothic{g}}
\newcommand{\gothh}{\gothic{h}}
\newcommand{\gotk}{\gothic{k}}
\newcommand{\gotm}{\gothic{m}}
\newcommand{\gotn}{\gothic{n}}
\newcommand{\gotp}{\gothic{p}}
\newcommand{\gotq}{\gothic{q}}
\newcommand{\gotr}{\gothic{r}}
\newcommand{\gots}{\gothic{s}}
\newcommand{\gotu}{\gothic{u}}
\newcommand{\gotv}{\gothic{v}}
\newcommand{\gotw}{\gothic{w}}
\newcommand{\gotz}{\gothic{z}}
\newcommand{\gotA}{\gothic{A}}
\newcommand{\gotB}{\gothic{B}}
\newcommand{\gotG}{\gothic{G}}
\newcommand{\gotL}{\gothic{L}}
\newcommand{\gotS}{\gothic{S}}
\newcommand{\gotT}{\gothic{T}}

\newcommand{\mn}{\marginpar{\hspace{1cm}*} }
\newcommand{\mnn}{\marginpar{\hspace{1cm}**} }

\newcommand{\mnq}{\marginpar{\hspace{1cm}*???} }
\newcommand{\mnnq}{\marginpar{\hspace{1cm}**???} }

\newcounter{teller}
\renewcommand{\theteller}{\Roman{teller}}
\newenvironment{tabel}{\begin{list}%
{\rm \bf \Roman{teller}.\hfill}{\usecounter{teller} \leftmargin=1.1cm
\labelwidth=1.1cm \labelsep=0cm \parsep=0cm}
                      }{\end{list}}

\newcounter{tellerr}
\renewcommand{\thetellerr}{(\roman{tellerr})}
\newenvironment{subtabel}{\begin{list}%
{\rm  (\roman{tellerr})\hfill}{\usecounter{tellerr} \leftmargin=1.1cm
\labelwidth=1.1cm \labelsep=0cm \parsep=0cm}
                         }{\end{list}}
\newenvironment{ssubtabel}{\begin{list}%
{\rm  (\roman{tellerr})\hfill}{\usecounter{tellerr} \leftmargin=1.1cm
\labelwidth=1.1cm \labelsep=0cm \parsep=0cm \topsep=1.5mm}
                         }{\end{list}}

\newcommand{\Ni}{{\bf N}}
\newcommand{\Ri}{{\bf R}}
\newcommand{\Ci}{{\bf C}}
\newcommand{\Ti}{{\bf T}}
\newcommand{\Zi}{{\bf Z}}
\newcommand{\Fi}{{\bf F}}

\newcommand{\proof}{\mbox{\bf Proof} \hspace{5pt}} 
\newcommand{\remark}{\mbox{\bf Remark} \hspace{5pt}}
\newcommand{\ruimte}{\vskip10.0pt plus 4.0pt minus 6.0pt}

\newcommand{\simh}{{\stackrel{{\rm cap}}{\sim}}}
\newcommand{\ad}{{\mathop{\rm ad}}}
\newcommand{\Ad}{{\mathop{\rm Ad}}}
\newcommand{\Aut}{\mathop{\rm Aut}}
\newcommand{\arccot}{\mathop{\rm arccot}}
\newcommand{\capp}{{\mathop{\rm cap}}}
\newcommand{\rcapp}{{\mathop{\rm rcap}}}
\newcommand{\diam}{\mathop{\rm diam}}
\newcommand{\divv}{\mathop{\rm div}}
\newcommand{\codim}{\mathop{\rm codim}}
\newcommand{\RRe}{\mathop{\rm Re}}
\newcommand{\IIm}{\mathop{\rm Im}}
\newcommand{\Tr}{{\mathop{\rm Tr}}}
\newcommand{\Vol}{{\mathop{\rm Vol}}}
\newcommand{\card}{{\mathop{\rm card}}}
\newcommand{\supp}{\mathop{\rm supp}}
\newcommand{\sgn}{\mathop{\rm sgn}}
\newcommand{\essinf}{\mathop{\rm ess\,inf}}
\newcommand{\esssup}{\mathop{\rm ess\,sup}}
\newcommand{\Int}{\mathop{\rm Int}}
\newcommand{\Leibniz}{\mathop{\rm Leibniz}}
\newcommand{\lcm}{\mathop{\rm lcm}}
\newcommand{\loc}{{\rm loc}}

\newcommand{\mod}{\mathop{\rm mod}}
\newcommand{\spann}{\mathop{\rm span}}
\newcommand{\one}{1\hspace{-4.5pt}1}

\newcommand{\DWR}{}

\hyphenation{groups}
\hyphenation{unitary}

\newcommand{\tfrac}[2]{{\textstyle \frac{#1}{#2}}}

\newcommand{\cb}{{\cal B}}
\newcommand{\cc}{{\cal C}}
\newcommand{\cd}{{\cal D}}
\newcommand{\ce}{{\cal E}}
\newcommand{\cf}{{\cal F}}
\newcommand{\ch}{{\cal H}}
\newcommand{\ci}{{\cal I}}
\newcommand{\ck}{{\cal K}}
\newcommand{\cl}{{\cal L}}
\newcommand{\cm}{{\cal M}}
\newcommand{\co}{{\cal O}}
\newcommand{\cs}{{\cal S}}
\newcommand{\ct}{{\cal T}}
\newcommand{\cu}{{\cal U}}
\newcommand{\cv}{{\cal V}}
\newcommand{\cx}{{\cal X}}
\newcommand{\cy}{{\cal Y}}
\newcommand{\cz}{{\cal Z}}

\newcommand{\wtozp}{W^{1,2}\raisebox{10pt}[0pt][0pt]{\makebox[0pt]{\hspace{-34pt}$\scriptstyle\circ$}}}
\newlength{\hightcharacter}
\newlength{\widthcharacter}
\newcommand{\covsup}[1]{\settowidth{\widthcharacter}{$#1$}\addtolength{\widthcharacter}{-0.15em}\settoheight{\hightcharacter}{$#1$}\addtolength{\hightcharacter}{0.1ex}#1\raisebox{\hightcharacter}[0pt][0pt]{\makebox[0pt]{\hspace{-\widthcharacter}$\scriptstyle\circ$}}}
\newcommand{\cov}[1]{\settowidth{\widthcharacter}{$#1$}\addtolength{\widthcharacter}{-0.15em}\settoheight{\hightcharacter}{$#1$}\addtolength{\hightcharacter}{0.1ex}#1\raisebox{\hightcharacter}{\makebox[0pt]{\hspace{-\widthcharacter}$\scriptstyle\circ$}}}
\newcommand{\scov}[1]{\settowidth{\widthcharacter}{$#1$}\addtolength{\widthcharacter}{-0.15em}\settoheight{\hightcharacter}{$#1$}\addtolength{\hightcharacter}{0.1ex}#1\raisebox{0.7\hightcharacter}{\makebox[0pt]{\hspace{-\widthcharacter}$\scriptstyle\circ$}}}

 \thispagestyle{empty}

\begin{center}
{\Large{\bf Degenerate elliptic operators  }}\\[2mm] 
{\Large{\bf in one dimension }}  \\[2mm]
\large  Derek W. Robinson$^1$ and Adam Sikora$^2$\\[2mm]

\normalsize{July 2009}
\end{center}

\vspace{5mm}

\begin{center}
{\bf Abstract}
\end{center}

\begin{list}{}{\leftmargin=1.8cm \rightmargin=1.8cm \listparindent=10mm 
   \parsep=0pt}
   \item
Let $H$ be the symmetric second-order differential operator on $L_2(\Ri)$ with domain $C_c^\infty(\Ri)$  and action  $H\varphi=-(c\,\varphi')'$
where  $ c\in W^{1,2}_{\rm loc}(\Ri)$ is a real function which is strictly positive  on $\Ri\backslash\{0\}$ but 
with $c(0)=0$.
We give a complete characterization of the self-adjoint extensions and the submarkovian extensions
of $H$.
In particular if $\nu=\nu_+\vee\nu_-$ where 
$\nu_\pm(x)=\pm\int^{\pm 1}_{\pm x} c^{-1}$ then  $H$ 
has a unique self-adjoint extension if and only if  $\nu\not\in L_2(0,1)$
and a unique submarkovian extension if  and only if $\nu\not\in L_\infty(0,1)$. 
In both cases the corresponding semigroup leaves $L_2(0,\infty)$ and $L_2(-\infty,0)$ invariant.

In addition we prove that for a general non-negative  $ c\in W^{1,\infty}_{\rm loc}(\Ri)$ the corresponding operator  $H$ has a unique submarkovian extension.
\end{list}

\vspace{5cm}

\noindent AMS Subject Classification: 47B25, 47D06, 35J15,  31C25.

\vspace{1cm}

\noindent
{\bf Home institutions:}    \\[3mm]
\begin{tabular}{@{}cl@{\hspace{10mm}}cl}
1. & Centre for Mathematics & 
  2. & Department of Mathematics  \\
&\hspace{15mm} and its Applications  & 
  &Macquarie University  \\
& Mathematical Sciences Institute & 
  & Sydney, NSW 2109  \\
& Australian National University& 
  & Australia \\
& Canberra, ACT 0200  & {}
  & \\
& Australia & {}
  & \\
& derek.robinson@anu.edu.au & {}
  &sikora@ics.mq.edu.au \\

\end{tabular}

\newpage

\setcounter{page}{1}

\null

\vspace{2mm}

\section{Introduction}\label{S1}
The theory of degenerate elliptic operators on $\Ri^d$  displays a number of significant features
which distinguishes it from the well understood non-degenerate theory.
If the degeneracies are weak then there is no great difference and the degenerate
theory can still be described by the techniques of the non-degenerate case, e.g.\
Harnack--Sobolev--Poincar{\'e} inequalities (see, for example,  \cite{Tru2} \cite{FKS}  \cite{FP}  \cite{Fra}   \cite{BM}  
  \cite{FLW}  \cite{SawW} and references therein).
If, however, the degeneracies are sufficiently strong the associated diffusion process
can exhibit  non-ergodic behaviour;  the degeneracies  can spontaneously introduce 
barriers and obstacles to the diffusion
\cite{ERSZ1} \cite{RSi}.
These properties can lead to quite unexpected phenomena such as cloaking and invisibility
 \cite{PSS} \cite{Wed}.

Despite the vast literature devoted to the subject 
many basic aspects of the degenerate theory are neither well developed nor well understood.  
For example, $ C_c^\infty(\Ri^d)$  is not necessarily a core for the degenerate operator acting on $L_2(\Ri^d)$ 
and one has  to consider boundary conditions  at the barriers and obstacles.
Moreover, these  boundary conditions can have a substantially different nature to the classical conditions of Dirichlet, Neumann or Robin. 
Therefore in this paper we give a detailed analysis of the simplest situation, divergence-form operators in one-dimension.
Although some  of the most interesting features are not apparent in one-dimension the analysis
does give a guide to possible features of the multi-dimensional case.

Let $x\in\Ri\mapsto c(x)$ be a real  function
in $ W^{1,2}_{\rm loc}(\Ri)$
which is strictly positive on $\Ri\backslash\{0\}$.
Define the second-order operator $H$ on $L_2(\Ri)$ with domain $D(H)=C_c^\infty(\Ri)$ by
\begin{equation}
H\varphi=-(c\,\varphi')'=-c\,\varphi''-c'\,\varphi'
\;.
\label{e4.2}
\end{equation}
Then $H$ is a positive-definite, symmetric, operator on $L_2(\Ri)$ with range in $L_1(\Ri)\cap L_2(\Ri)$.
Our aim is to study the self-adjoint  extensions  of $H$ and the semigroups they generate.
In particular we are interested in the submarkovian extensions, i.e.\ the extensions which generate 
submarkovian semigroups.
One such extension always exists because  the
 closure $\overline h$ of the quadratic form $h$ associated with $H$,
i.e.\ the form
\begin{equation}
h(\varphi)=(\varphi,H\varphi)=\int_\Ri dx\,c(x)\,|\varphi'(x)|^2
\label{e4.2.1}
\end{equation}
with domain $D(h)=C_c^\infty(\Ri)$, is a Dirichlet form.
Therefore the corresponding self-adjoint extension $H_F$ of $H$, the Friedrichs extension,
is submarkovian. 
(For  background on submarkovian semigroups and Dirichlet forms see \cite{FOT} \cite{BH} \cite{MR}.)

If   $c>0$ on the whole line then $H$ is essentially self-adjoint (see \cite{DunS2}, Corollary ~XIII.6.15) but the 
situation is complicated by a degeneracy at the origin.
Then there is a trichotomy of self-adjoint extensions which can be indexed by  the local properties of the 
functions
  $\nu_\pm$   defined by
\begin{equation}
x>0\mapsto \nu_+(x)=
\int^1_x ds\,c(s)^{-1}
\;\;\;\;\;\;
{\rm and }\;\;\; \;\;\;x>0\mapsto \nu_-(x)=
\int^{-x}_{-1}ds\,c(s)^{-1}
\;.
\label{e4.3.0}
\end{equation}
These functions are $H$-harmonic, i.e.\ $H\nu_\pm=0$.
The choice $\pm\,1$ as integral limits in their definition is arbitrary since  only the behaviour of the functions
at the origin is important.

\begin{thm}\label{text4.1}
$\;\;$Assume $c(0)=0$.
\begin{tabel}
\item\label{text4.1-1}
If $\nu_+\vee\nu_-\not\in L_2(0,1)$     then  $H$ is essentially self-adjoint
 and the self-adjoint  closure $\overline H$ of $H$ generates a  submarkovian semigroup $S$
 which leaves $L_2(-\infty,0)$ and $L_2(0,\infty)$ invariant.
\item\label{text4.1-2}
If  $\nu_+\vee\nu_-\in L_2(0,1)$ but  $\nu_+\vee\nu_-\not\in  L_\infty(0,1)$ 
then $H$  has a one-parameter family of self-adjoint extensions
but only one extension
generates a positive semigroup.
This semigroup is submarkovian and 
leaves  $L_2(-\infty,0)$ and $L_2(0,\infty)$ invariant.
\item\label{text4.1-3} If $\nu_+\vee\nu_-\in  L_\infty(0,1)$  then $H$  has a one-parameter family of self-adjoint extensions
none of which leave $L_2(-\infty,0)$ and $L_2(0,\infty)$ invariant.
Moreover, there is a one-parameter subfamily of submarkovian extensions.
\end{tabel}

In   particular $H$ has a unique self-adjoint extension if and only if $\nu_+\vee\nu_-\not\in L_2(0,1)$  and a unique submarkovian extension if and only if
$\nu_+\vee\nu_-\not\in L_\infty(0,1)$.
\end{thm}

The self-adjoint extensions of $H$ can in part be described by classical boundary conditions.
In Case~\ref{text4.1-1} all elements in  the domain of the unique self-adjoint extension satisfy the condition $(c\,\varphi')(0_\pm)=0$.
In Case~\ref{text4.1-2} the same conditions characterize the domain of the unique submarkovian extension.
In Case~\ref{text4.1-3}  the boundary condition for a general self-adjoint extension is given by
\[
\beta\,\Big((c\,\varphi')(0_+)-(c\,\varphi')(0_-)\Big)=\alpha\,\Big(\varphi(0_+)-\varphi(0_-)\Big)
\]
where $\alpha,\beta\in \Ri^2\backslash(0,0)$  and   
the  submarkovian extensions are determined by the condition $\alpha\, \beta\geq0$.
 The non submarkovian extensions in Case~\ref{text4.1-2} are exceptional.
There is no  comparable classification of these extensions
and they  do not have any obvious probabilistic interpretation.

The uniqueness criteria in the last statement of Theorem~\ref{text4.1}  are a measure of the order of degeneracy 
of $c$ at the origin.
If $c(x)=O(x^{\delta_\pm})$ as $x\to0_\pm$ then there is a unique
selfadjoint  extension if and only if either $\delta_+\geq 3/2$ or $\delta_-\geq 3/2$
and a unique submarkovian extension if and only if either $\delta_+\geq 1$ or $\delta_-\geq 1$.
It is notable that it suffices to have  a `strong' degeneracy on one side.
This one-sideness has been stressed by Weder in the context of cloaking \cite{Wed}, Theorem~2.5.
Note also that these criteria are local properties and do not  depend on the behaviour of $c$ at infinity.
Next we examine a different type of characterization  of uniqueness.

It follows from general operator theory that $H$  has a unique self-adjoint
extension, if and only if  the range of $(I+H)$ is dense in $L_2(\Ri)$.
There is a similar characterization of  uniqueness of the submarkovian extension by an
 $L_1$-range condition at least if   the  coefficient $c$ satisfies a growth estimate.
Define   the positive increasing functions
\begin{equation}
x\geq 1\mapsto \mu_+(x)=\int^x_1ds\,s\,c(s)^{-1}\;\;\;\;\;{\rm and}\;\;\;\;\;x\geq 1\mapsto \mu_-(x)=-\int^{-1}_{-x}ds\,s\,c(s)^{-1}
\;.
\label{e4.3.00}
\end{equation}
Then $H\mu_\pm=1$.
Since  $H$  is a second-order elliptic operator it is  both dissipative and dispersive
as an operator on $L_1(\Ri)$.
(see Section~\ref{S5}).
In particular it is $L_1$-closable.
Then the closure generates a positive contractive
semigroup on $L_1(\Ri)$ if and only if the range of $(I+H)$ is $L_1$-dense.

\begin{thm}\label{cext1.0}
Consider the following conditions.
\begin{tabel}
\item\label{cext1.0-3}
$(I+H)C_c^\infty(\Ri\backslash\{0\})$ is dense in $L_1(\Ri)$.
\item\label{cext1.0-2}
$(I+H)C_c^\infty(\Ri)$ is dense in $L_1(\Ri)$.
\item\label{cext1.0-1}
The operator $H$ has a unique submarkovian extension.
 \end{tabel}
 Then  {\rm\ref{cext1.0-3}$\Rightarrow$\ref{cext1.0-2}$\Rightarrow$\ref{cext1.0-1}} and if $\mu_+\wedge \mu_-\not\in L_\infty(1,\infty)$
then  {\rm \ref{cext1.0-1}$\Rightarrow$\ref{cext1.0-3}}.
\end{thm}

The proof of Theorem~\ref{text4.1} is  in two steps. 
First, in Section~\ref{S2},   we consider the analogous problems  on the left and right half-lines.
Secondly,  in Section~\ref{S3},  we marry  together  the results
for the two half-lines to obtain the description of the various extensions   on the line.
 Our analysis on the half-line  overlaps with the early work of Feller \cite{Feller} \cite{Feller2} \cite{Feller3}
(see also \cite{Mandl})  
but our emphasis  is  different and the arguments  are independent.
Feller classified extensions of operators acting  on $L_1$, or on  $C_b$, 
which generate positive contraction semigroups. 
The boundary conditions $\nu_\pm\not\in L_\infty(0,1)$ are interpretable in Feller's terminology.
The condition $\nu_+\not\in L_\infty(0,1)$   states that $0_+$  is  an inaccessible natural boundary
 (see \cite{Mandl}, pages~24--25).

The $L_2$-theory has, however,  several different features not shared by the  $L_1$-theory  since  there are $L_2$-extensions which generate continuous 
 semigroups which do not extend to  $L_1$ or  $L_\infty$.
 Moreover, the $L_2$-arguments do not require any growth restrictions on $c\,$; there are no boundary conditions at infinity.
 The proof of Theorem~\ref{cext1.0}, which is given in 
 Section~\ref{S5}, is,  however,  based on $L_1$-arguments which depend in part on growth properties.
The growth condition, $\mu_+\wedge \mu_-\not\in L_\infty(1,\infty)$,    coincides with  Feller's criterion for 
$\pm\,\infty$  to  be inacessible boundaries.
The range conditions \ref{cext1.0-3} and \ref{cext1.0-2} in Theorem~\ref{cext1.0} are of independent interest as they imply that the submarkovian 
semigroup is conservative.

Our arguments extend  to operators defined on finite intervals which are degenerate at both endpoints  \cite{Ulm} \cite{CMP}.
This is briefly discussed in   Section~\ref{S4}
where we  establish the following simple statement  for operators with Lipschitz continuous coefficients
for which the zero set might be quite complicated.

\begin{thm}\label{cext1.1} If $c\in W_{\rm loc}^{1,\infty}(\Ri)$ is non-negative
then $H$ has a unique submarkovian extension.
\end{thm}

The proof of Theorem~\ref{cext1.1} uses a mixture of $L_1$- and $L_2$-arguments.
But these are all of a local nature and 
again no growth condition at infinity is necessary.

Theorems~\ref{cext1.0} and \ref{cext1.1}  should have analogues in higher dimensions.

\section{The half-line}\label{S2}

In  this section we examine the self-adjoint extensions of the restriction $H_+=H|_{C_c^\infty(0,\infty)}$
of $H$ to the right half-line. 
The analysis of $H_-=H|_{C_c^\infty(-\infty,0)}$ is similar.

First,  since $c>0$ on $\Ri\backslash\{0\}$ the domain of the adjoint $H_+^*$ of  $H_+$  is given by 
\begin{equation}
D(H_+^*)=\{\varphi\in L_2(0,\infty)\cap AC_{\rm loc}(0,\infty): c\,\varphi'\in AC_{\rm loc}(0,\infty)\,,\,( c\, \varphi')'\in L_2(0,\infty)\}
\label{esm2.1}
\end{equation}
and $H_+^*\varphi=-(c\,\varphi')'$ for $\varphi\in D(H_+^*)$.
(The details of this identification are given  in \cite{Kat1}, Sections~III.2.3 and III.5.5, under slightly stronger assumptions on $c$.)
Secondly, the domain of the closure $\overline H_+$ of $H_+$  is  obtained as the restriction of  $D(H_+^*)$ by a boundary condition
(see, for example, \cite{Sto}, Theorem~10.11, \cite{DunS2}, Section~XIII.2 or \cite{Far}, Section~13). 
In principle the boundary condition is the direct sum of a boundary value at the origin and a boundary value at infinity.
But by an argument of Wintner (see  \cite{DunS2}, Theorem~XIII.6.14) there is no boundary value at infinity.
Therefore 
\begin{equation}
D(\overline H_+)=\{\varphi\in D(H_+^*): B_+(\varphi,\psi)=0 \mbox{ for all } \psi\in D(H_+^*)\}
\;,
\label{ebc}
\end{equation}
where the boundary value  $B_+(\,\cdot\,,\,\cdot\,)$ is a bilinear functional over $D(H_+^*)$ defined by
\begin{eqnarray*}
B_+(\varphi,\psi)=(H_+^*\varphi,\psi)-(\varphi,H_+^*\psi)
=\lim_{x\to0_+}\Big(({c\,\varphi'})(x)\,\psi(x)-\varphi(x)\,{ (c\,\psi')}(x)\Big)
\end{eqnarray*}
for all $\varphi,\psi\in D(H_+^*)$.
The limit in the boundary term exists although the limits of the individual terms in the expression
do not necessarily exist.
In the sequel we  write $\chi(0_+)=\lim_{x\to0_+}\chi(x)$ for 
any function $\chi$  over  $\langle0,\infty\rangle$.
It is implicit in this usage that the limit  does exist.

\begin{prop}\label{pext2.1}$\;$Let $\nu_+$ be the harmonic function defined on $\langle0,\infty\rangle$ by $(\ref{e4.3.0})$.
\begin{tabel}
\item\label{pext2.1-1}
If $\nu_+\not\in L_2(0,1)$ then 
$H_+$ is essentially self-adjoint.
\item\label{pext2.1-2}
If $\nu_+\in L_2(0,1)$ then 
\begin{equation}
D(\overline H_+)=\{\varphi\in D(H_+^*):  
\varphi(0_+)=0=(\nu_+\,c\,\varphi')(0_+)\}\subsetneqq D(H_+^*)
\label{eext2.2}
\end{equation}
and $H_+$ has deficiency indices $(1,1)$.
\end{tabel}
\end{prop}
\noindent{\bf Proof of Proposition~\ref{pext2.1}.\ref{pext2.1-1}}$\;$
The proof relies on the following lemma.

\begin{lemma}\label{lext2.1}
Assume $\nu_+\not\in L_2(0,1)$.
If $\varphi\in D(H_+^*)$ then 
$(c\,\varphi')(0_+)=0=(c\,\varphi\,\varphi')(0_+)$.
\end{lemma}
\proof\
If $\varphi\in D(H_+^*)$ then   
$\lim_{x\to0_+}(c\, \varphi')(x)=\varepsilon $ exists
and  $\varphi'\sim \varepsilon\,c^{-1}$ on $\langle0,1]$.
But $\varphi\not\in L_2(0,1)$ unless $\varepsilon=0$.
So  $(c\,\varphi')(0_+)=0$,

\[
(c\,\varphi')(x)=\int^x_0ds\,(c\,\varphi')'(s)=-\int^x_0ds\,(H_+^*\varphi)(s)
\]
and  $(c\,\varphi')(x)=O(x^{1/2})$ as $x\to0_+$.
Next we argue that $\lim_{x\to0_+}(c\,\varphi\,\varphi')(x)$ exists.

The coefficient $c$ is strictly positive and bounded on each bounded  interval $[a,b\,]$ with $b>a>0$.
Since $\varphi\in D(H^*)$
one  has $\varphi'\in L_2(a,b)$.
Hence
\begin{equation}
\int^b_a(H_+^*\varphi)\,\varphi-\int^b_a c\,|\varphi'|^2=(c\,\varphi\,\varphi')(a)-(c\,\varphi\,\varphi')(b)
\;.
\label{eext2.20}
\end{equation}
But the limit as $a\to0_+$ of the left hand side exists although it is not necessarily finite. 
Therefore $(c\,\varphi\,\varphi')(0_+)$ exists
and we next  argue that it is zero.

Suppose $(c\,\varphi\,\varphi')(0_+)\geq \varepsilon>0$.
Thus there is a $\delta>0$ such that $(c\,\varphi\,\varphi')(x)\geq \varepsilon/2$ for $x\in\langle0,\delta\,]$.
Moreover, $(c\,\varphi')(x)=O(x^{1/2})$ as $x\to0_+$ by the foregoing argument.
Therefore $|\varphi(x)|\geq a\,x^{-1/2}$ as $x\to0_+$ with $a>0$ unless $\varepsilon=0$.
But then $\varepsilon$ must be zero since  $\varphi\in L_2(0,1)$.
An identical argument applies if $(c\,\varphi\,\varphi')(0_+)$ is initially assumed to be negative.
\hfill$\Box$

\bigskip

The symmetric operator  $H_+$ is essentially self-adjoint if and only if  $(I+H_+)C_c^\infty(0,\infty)$ is dense in $L_2(0,\infty)$.
Assume there is a $\psi\in L_2(0,\infty)$ such that $(\psi,(I+H_+)\varphi)=0$ for all $\varphi\in C_c^\infty(0,\infty)$.
Then $|(\psi,H_+\varphi)|\leq \|\psi\|_2\,\|\varphi\|_2$ so $\psi\in D(H_+^*)$ and $(I+H_+^*)\psi=0$.
But 
\[
\int^b_0(H_+^*\psi)\psi=-(c\,\psi\,\psi')(b)+\int^b_0 c\,|\psi'|^2\geq -(c\,\psi\,\psi')(b)=-2^{-1}(c\,(\psi^2)')(b)
\]
because $(c\,\psi\,\psi')(0_+)=0$ by Lemma~\ref{lext2.1}.
Since $\psi$ is square-integrable  $\psi^2$ cannot be monotone increasing.
Hence the derivative of $\psi^2$ must take non-positive values for large $b$, i.e.\ there is a sequence
$b_n\to\infty$ for which $-(c\,(\psi^2)')(b_n)\geq0$.
Then in the limit $n\to\infty$ one concludes that $(H_+^*\psi,\psi)\geq0$.
Therefore
\[
\|\psi\|_2^2\leq ((I+H_+^*)\psi,\psi)=0
\]
and so $\psi=0$.
Thus the range of $(I+H_+)$ is dense, $\overline H_+$ is self-adjoint and $D(\overline H_+)=D(H_+^*)$.
This completes the proof of  Proposition~\ref{pext2.1}.\ref{pext2.1-1}.\hfill$\Box$

\bigskip

The proof of  Proposition~\ref{pext2.1}.\ref{pext2.1-2} is based on the following
three  lemmas.
The first  does not require any special assumption on the behaviour of $\nu_+$.

\begin{lemma}\label{lext2.3}
If $\varphi\in D(\overline H_+)$ then $|(c\,\varphi')(x)|\leq x^{1/2}\,\|\overline H_+\varphi\|_2$ for all $x\in\langle0,1]$.
In particular  $(c\,\varphi')(0_+)=0$.
\end{lemma}
\proof\
Fix $\varphi\in D(\overline H_+)$.
Choose a sequence $\varphi_n\in C_c^\infty(0,\infty)$ such that $\|\varphi_n-\varphi\|_2\to0$  and $\|H_+\varphi_n'-\overline H_+\varphi\|_2\to0$ as $n\to\infty$.
Then $(c\,\varphi_n')(x)=-\int^x_0\,H_+\varphi_n$ and 
\[
|(c\,\varphi_n')(x)-(c\,\varphi_m')(x)|\leq x^{1/2}\|H_+(\varphi_n-\varphi_m)\|_2
\;.
\]
Thus $c\,\varphi_n'$ converges uniformly on $\langle0,1]$  to a limit $\psi$.
But $c$ is strictly positive on each closed interval $I\subset \langle0,1]$.
Hence $\varphi_n'$ converges uniformly on $I$ to $c^{-1}\psi$.
In particular it is $L_2(I)$-convergent to $c^{-1}\psi$.
Since $\varphi_n$ is $L_2(I)$-convergent to $\varphi$ and the maximal operator of differentiation
is closed on $L_2(I)$ it follows that $\psi=c\,\varphi'$.
Therefore 
\begin{equation}
(c\,\varphi')(x)=\lim_{n\to\infty}(c\,\varphi_n')(x)=-\lim_{n\to\infty}\int^x_0ds\,( H_+\varphi_n)(s)=-\int^x_0ds\,(\overline H_+\varphi)(s)\;.
\label{e3.2.1}
\end{equation}
Hence $|(c\,\varphi')(x)|\leq x^{1/2}\,\|\overline H_+\varphi\|_2$.
\hfill$\Box$

\bigskip

The second lemma gives control over the singularity of $\nu_+$.

\begin{lemma}\label{lext2.11}
Assume $\nu_+\in L_2(0,1)$.
Then $x^{1/2}\nu_+(x)\to0$ as $x\to0_+$.
\end{lemma}
\proof\
Suppose the statement is false.
Then there exists a decreasing sequence $1\geq x_1\geq x_2\ldots>0$  and an
$\varepsilon>0$ such that 
$x_n^{1/2} \,\nu_+(x_n) \geq \varepsilon>0$ for all $ n\geq 1$.
But by passing to a  subsequence if necessary one may assume that 
$x_n > 2\,x_{n+1}$. 
Since $\nu_+(x)$ decreases with $x$ one then has
\[
\int^1_0ds\, \nu_+(s)^2 \geq \sum_{n\geq1} \nu(x_n)^2(x_{n}-x_{n+1}) \geq  \sum_{n\geq1} \nu_+(x_n)^2x_{n}/2  \geq  (\varepsilon^2/2)\sum_{n\geq1}\;1
\]
which contradicts the assumption $\nu_+\in L_2(0,1)$.
\hfill$\Box$

\bigskip

Combination of the foregoing lemmas leads to the following conclusion.

\begin{cor}\label{cext2.1}
Assume $\nu_+\in L_2(0,1)$.
If  $\varphi\in D(\overline H_+)$ then $(\nu_+\,c\,\varphi')(0_+)=0$.
\end{cor}

The final preparatory lemma is the following.

\begin{lemma}\label{lext2.2}
Assume  $\nu_+\in L_2(0,1)$. 
If  $\varphi\in D(\overline H_+)$ then $\varphi(0_+)=0$.
\end{lemma}
\proof\ 
Fix  $\varphi\in D(\overline H_+)$.
Again there is a sequence $\varphi_n\in C_c^\infty(0,\infty)$ such that $\|\varphi_n-\varphi\|_2\to0$  and $\|H_+\varphi_n'-\overline H_+\varphi\|_2\to0$  as $n\to\infty$.
But
\[
\varphi_n'(x)=c(x)^{-1}\int^x_0ds\,(c\,\varphi_n')'(s)
\]
for all $x>0$.
Therefore 
\begin{equation}
\varphi_n(x)=\int^x_0ds\,c(s)^{-1}\int^s_0dt\, (c\,\varphi_n')'(t)=-\int^x_0dt\,\nu_x(t)(H_+\varphi_n)(t)
\label{esm3.1.1}
\end{equation}
where $\nu_x(t)=\int^x_tds\,c(s)^{-1}$.
But $0\leq \nu_x(t)\leq \nu_+(t)$ for all $x,t\in\langle0,1]$
so
\begin{equation}
|\varphi_n(x)|\leq \|\nu_+\|_2\|H_+\varphi_n\|_2
\label{esm3.9}
\end{equation}
for all $x\in\langle0,1]$.
Similarly
\begin{equation}
|\varphi_n(x)-\varphi_m(x)|\leq \|\nu_+\|_2\|H_+(\varphi_n-\varphi_m)\|_2
\label{esm3.91}
\end{equation}
for all $x\in\langle0,1]$.
It follows that $\varphi_n$ converges to $\varphi$ uniformly on compact subsets of $\langle0,1]$.
But then one deduces from (\ref{esm3.1.1}) that
\begin{equation}
\varphi(x)=-\int^x_0dt\,\nu_x(t)(\overline H_+\varphi)(t)
\;.\label{esm3.10}
\end{equation}
Hence
\[
|\varphi(x)|\leq \|\nu_+\|_2 \bigg(\int^x_0dt\,|(\overline H_+\varphi)(t)|^2\bigg)^{1/2}
\]
for all $x\in\langle0,1]$.
Therefore $\lim_{x\to0_+}\varphi(x)=0$.\hfill$\Box$

\bigskip

\noindent{\bf Proof of Proposition~\ref{pext2.1}.\ref{pext2.1-2}}$\;$
Suppose $\nu_+\in L_2(0,1)$.
Then $D(\overline H_+)\subseteq D_0$ where
\[
D_0=\{\varphi\in D(H_+^*): \varphi(0_+)=0=(\nu_+\,c\,\varphi')(0_+)\}
\]
by Lemma~\ref{lext2.2} and Corollary~\ref{cext2.1}.
Next we prove the converse inclusion.

 Fix $\varphi\in D_0$. 
If $\psi\in D(H_+^*)$ then  
$\lim_{x\to0_+}(c\, \psi')(x)$ exists.
Therefore there is a $b>0$ such that $|(c\, \psi')(x)|\leq b$ and $|\psi(x)|\leq  b\,\nu_+(x)$
for  all small $x$.
But then 
\[
|B_+(\varphi,\psi)|\leq b\,|\varphi(0_+)|+b\,|(\nu_+\,c\,\varphi')(0_+)|=    0
\]
and  so  $\varphi\in D(\overline H_+)$.
Therefore $D_0=D(\overline H_+)$.
Then it follows from (\ref{esm2.1}) that $D_0$ is a strict subset of $D(H_+^*)$.
Hence $D(\overline H_+)$ is a strict subset of $D(H_+^*)$ and $H_+$ must have deficiency indices $(1,1)$.
\hfill$\Box$

\bigskip

At this stage we can prove an analogue of Theorem~\ref{text4.1}
for operators on the half-line.

\begin{thm}\label{tnsm1.1}$\;\;\;$
\begin{tabel}
\item\label{tnsm1.1-1}
If $\nu_+\not\in L_2(0,1)$  then  $H_+$ is essentially self-adjoint
 and its   closure   generates a  submarkovian semigroup.
\item\label{tnsm1.1-2}
If $\nu_+\in L_2(0,1)$ but $\nu_+\not\in L_\infty(0,1)$
then $H_+$ has  a one-parameter family of self-adjoint extensions
each of which  generates a positive semigroup
but  only one extension, corresponding  to the boundary condition $(c\,\varphi')(0_+)=0$,
 generates a submarkovian semigroup.
\item\label{tnsm1.1-3}
If $\nu_+\in  L_\infty(0,1)$ then   $H_+$ has a one-parameter family of self-adjoint extensions
 characterized by the  classical Dirichlet,  Neuman and Robin boundary conditions.  
All the extensions generate positive semigroups
 and the  positive$($-definite$)$ extensions generate submarkovian semigroups.
\end{tabel}

In particular $H_+$ has a unique self-adjoint extension if and only if $\nu_+\not\in L_2(0,1)$ and a unique submarkovian extension if and only if $\nu_+\not\in L_\infty(0,1)$.
\end{thm}
\proof\
Theorem~\ref{tnsm1.1}.\ref {tnsm1.1-1} is a direct consequence of the first statement of Proposition~\ref{pext2.1}.
Since $H_+$ is essentially self-adjoint the self-adjoint closure must coincide with the Friedrichs extension which generates
a submarkovian semigroup as remarked in Section~\ref{S1}.

The second and third statements of the theorem require more detailed analysis of the self-adjoint extensions of $H_+$.
Throughout the following we assume $\nu_+\in L_2(0,1)$.

\smallskip

It follows from Proposition~\ref{pext2.1}.\ref{pext2.1-2} that the deficiency indices of $H_+$ are $(1,1)$.
Therefore  the codimension of $D(\overline H_+)$  in $D(H^*)$ is two
and  $D(H_+^*)$ can be spanned by $D(\overline H_+)$ and two auxiliary functions
which we choose to be local solutions of the harmonic equation $H_+^*\psi=0$.

Let   $\sigma_+\in C_c^\infty(0,\infty)$ satisfy  $\sigma_+(x)=1$ if $x\in[0,1\rangle$ and $\sigma_+(x)=0$ if $x\geq 2$.
(The choice of values $1$ and $2$ is not significant. 
One could equally well assume that $\sigma_+(x)=1$ if $x\in[0,\varepsilon\rangle$ and $\sigma_+(x)=0$ if $x\geq \delta$
with $0<\varepsilon<\delta$.
It is only important that $\sigma_+$ is equal to one near the origin.)
Next   set $\tau_+=\nu_+\,\sigma_+$.
Clearly one has  $\sigma_+,\tau_+\in D(H^*)$ and $(H_+^*\sigma_+)(x)=0=(H_+^*\tau_+)(x)$ for $x\in\langle 0,1\rangle$.
But it follows from (\ref{eext2.2}) that  $\sigma_+,\tau_+\not\in D(\overline H_+)$.
In the sequel $(\sigma_+,\tau_+)$ always denotes a pair of functions constructed in this manner.

\begin{prop}\label{p3.2.1}
Assume $\nu_+\in L_2(0,1)$. 
Then 
$D(H_+^*)=D(\overline H_+)+\spann \sigma_+ +\spann \tau_+$.
\end{prop}
\proof\ 
First, by definition $D(\overline H_+)+\spann \sigma_+ +\spann \tau_+\subseteq D(H_+^*)$.

Secondly, fix $\Phi\in D(H_+^*)$.
Then 
$\lim_{x\to0_+}(c\,\Phi')(x)=b$ exists.
But if $\Psi=\Phi+b\,\tau_+$ then $\Psi\in D(H_+^*)$ and $(H_+^*\Psi)(x)=(H_+^*\Phi)(x)$
for $x\in\langle0,1\rangle$.
Moreover, $(c\,\Psi')(0_+)=(c\,\Phi')(0_+)-b=0$. 
Then, however,
\begin{eqnarray*}
\Psi(x)&=&\Psi(1)-\int^1_xds\,c(s)^{-1}\int^s_0dt\,(H_+^*\Psi)(t)\\[5pt]
&=&\Psi(1)-\int^1_0ds\,c(s)^{-1}\int^s_0dt\,(H_+^*\Phi)(t)+\int^x_0ds\,c(s)^{-1}\int^s_0dt\,(H_+^*\Phi)(t)
\end{eqnarray*}
for all $x\in\langle0,1\rangle$.
Both the latter integrals are well defined since $\nu_+\in L_2(0,1)$.
Moreover, \begin{eqnarray*}
\Big|\int^x_0ds\,c(s)^{-1}\int^s_0dt\,(H_+^*\Phi)(t)\Big|&=&\Big|\int^x_0dt\,(H_+^*\Phi)(t)\int^x_tc(s)^{-1}\Big|\\[5pt]
&\leq&\|\nu_+\|_2\bigg(\int^x_0dt\,|(H_+^*\Phi)(t)|^2\bigg)^{1/2}
\end{eqnarray*}
which establishes that \[
a=\lim_{x\to0_+}\Psi(x)=\lim_{x\to0_+}\Big(\Phi(x)+b\,\tau_+(x)\Big)
\]
exists and $a<\infty$.
Next define $\varphi$ by 
\[
\varphi=\Phi-a\,\sigma_++b\,\tau_+=\Psi-a\,\sigma_+
\;.
\]
Then $\varphi\in D(H_+^*)$.
Moreover, if $x\in\langle0,1\rangle$ then 
\[
(c\,\varphi')(x)=(c\,\Psi')(x)=-\int^x_0ds\,(H_+^*\Psi)(s)
\;.
\]
Therefore $|(c\,\varphi')(x)|=O(x^{1/2})$ as $x\to0_+$.
Hence $(\nu_+\,c\,\varphi')(0_+)=0$ by  Lemma~\ref{lext2.11}.
But $\varphi(0_+)=\Psi(0_+)-a=0$ by the definition of $a$.
Therefore $\varphi\in D(\overline H_+)$ by Proposition~\ref{pext2.1}.\ref{pext2.1-2}.
Hence $D(H_+^*)\subseteq D(\overline H_+)+\spann \sigma_+ +\spann \tau_+$.
\hfill$\Box$

\bigskip

The one-parameter family of self-adjoint extensions of $\overline H_+$  can  now be specified by restricting  $H_+^*$ to subspaces of $D(H_+^*)$ obtained by supplementing 
$D(\overline H_+)$ through the addition of a one-dimensional subspace of $\spann \sigma_++\spann\tau_+$.

\begin{prop}\label{p3.2.2}
Assume  $\nu_+\in L_2(0,1)$.
Let $(\alpha,\beta)\in\Ri^2\backslash(0,0)$.
Define 
\[
D_{\alpha,\beta}=D(\overline H_+)+\spann(\beta\,\sigma_+-\alpha\,\tau_+)=\{\varphi\in D(H^*): B_+(\beta\,\sigma_+-\alpha\,\tau_+, \varphi)=0\}
\;.
\]
Then the restriction $H_{\alpha,\beta}$ of $H_+^*$ to $D_{\alpha,\beta}$ is a  self-adjoint extension of  $\overline H_+$ and 
 $H_{\alpha,\beta}=H_{\alpha',\beta'}$ if  and only if $\alpha\,\beta'=\beta\,\alpha'$.
\end{prop}
\proof\
It follows by the definition of  $H_{\alpha,\beta}$ that  $\overline H_+\subseteq H_{\alpha,\beta}\subseteq H_+^*$. 
Therefore $\Phi\in D(H_{\alpha,\beta}^*)$ if and only if $\Phi\in D(H_+^*)$ and $B_+(\Psi,\Phi)=0$ for all $\Psi\in D_{\alpha,\beta}$
(see the lemma on page 86 of \cite{Far}).
But if $\Phi=\varphi+a\,\sigma_+-b\,\tau_+\in D(H_+^*)$   and $\Psi=\psi+\lambda\,(\beta\,\sigma_+-\alpha\,\tau_+)\in D_{\alpha,\beta}$
with  $\varphi,\psi\in D(\overline H_+)$ and $a,b,\lambda\in\Ri$ then
\[
B_+(\Psi,\Phi)=\lambda\,B_+(\beta\,\sigma_+-\alpha\,\tau_+, \Phi)=\lambda\,(a\,\alpha-b\,\beta)
\;.
\]
Thus $\Phi\in  D(H_{\alpha,\beta}^*)$ if and only if $a\,\alpha=b\,\beta$.
But then $\Phi\in D_{\alpha,\beta}$.
Since $H^*_{\alpha,\beta}\subseteq H^*$ it follows that $H^*_{\alpha,\beta}\Phi=H_+^*\Phi=H_{\alpha,\beta}\Phi$
and consequently  $H_{\alpha,\beta}$ is self-adjoint.

The definitions of $D_{\alpha,\beta}$ and $H_{\alpha,\beta}$  are clearly independent of a change $(\alpha,\beta)\to(\alpha',\beta')$ if 
$\alpha\,\beta'=\beta\,\alpha'$.
\hfill$\Box$

\bigskip

Next we prove the second statement of Theorem~\ref{tnsm1.1}.
This  is the most complicated and the most interesting case.
Throughout the sequel we use the notation  $h_+$ and $h_{\alpha,\beta}$ for the closed quadratic forms associated with the 
operators $H_+$ and $H_{\alpha,\beta}$.

\begin{thm}\label{c3.2.2}
Assume  $\nu_+\in L_2(0,1)$ but $\nu_+\not\in L_\infty(0,1)$. 
Then $\overline h_+=h_{0,1}$, the Friedrichs  extension $H_F$ of $H_+$ is equal to $H_{0,1}$ and
\[
D(H_F)=\{\varphi\in D(H_+^*): (c\, \varphi')(0_+)=0\}
\;.
\]

The semigroup generated by $H_F$ is submarkovian.
Moreover, $H_F$  is the unique self-adjoint extension of $H_+$ which generates a submarkovian semigroup.
\end{thm}

A key  part  of the proof is the following lemma which follows from the arguments 
of  \cite{ERSZ1}, Proposition~6.5, or   \cite{RSi}, Example~3.3.

\begin{lemma} \label{lext3.137}
Assume $\nu_+\not\in L_\infty(0,1)$.
Then $\sigma_+\in D(\overline h_+)$.
\end{lemma}
\proof\
The proof is a repetition of the arguments of \cite{ERSZ1} but since the result is crucial for
 the sequel we sketch the details.

Define $\chi_n \colon \Ri \to [0,1]$ by
\begin{equation}
\chi_n(x)
= \left\{ \begin{array}{ll}
   1 & \mbox{if } x\leq n^{-1} \, ,  \\[5pt]
   \nu_n^{-1} \, \nu_+(x) & \mbox{if } x\in\langle n^{-1},1\rangle \, , \\[5pt]
   0 & \mbox{if } x \geq 1 \, ,\label{apid}
          \end{array} \right.
\end{equation}
where $n\in\Ni$ and $\nu_n=\nu_+(n^{-1})$.
Since $\nu_+\not\in L_\infty(0,1)$ it follows that $\nu_n\to\infty$ as $n\to\infty$.
Now  $\chi_n$ is absolutely continuous,
 $\|\chi_n\,\sigma_+\|_2\to 0$ as $n\to\infty$ and
 $\chi_n\,\sigma_+'=0$.
 Moreover,
\[
\int^\infty_0 c\,|(\chi_n\,\sigma_+)'|^2=\int^1_{n^{-1}}c\,|\chi_n'\,\sigma_+|^2=\int_{n^{-1}}^1 c\,|\chi_n'|^2
\]
because $\sigma_+=1$ on $[0, 1\rangle$.
But  $\chi_n'=\nu_n^{-1}c^{-1}$ on the interval $[n^{-1},1]$.
Therefore
\[
\int^\infty_0 c\,|(\chi_n\,\sigma_+)'|^2=\nu_n^{-2}  \int_{n^{-1}}^1 c^{-1}= \nu_n^{-1}
\;.
\]
Next  $\supp(1-\chi_n)\,\sigma_+\subseteq [n^{-1},1]$.
Consequently one checks that  
$(1-\chi_n)\,\sigma_+\in D(\overline h_+)$.
Now $\|(1-\chi_n)\,\sigma_+-\sigma_+\|_2=\|\chi_n\,\sigma_+\|_2\to 0$ as $n\to\infty$
and 
\[
\overline h_+((1-\chi_n)\,\sigma_+-(1-\chi_m)\,\sigma_+)\leq 2\int^\infty_0 c\,|(\sigma_+\,\chi_n)'|^2+2\int^\infty_0 c\,|(\sigma_+\,\chi_m)'|^2=2\,(\nu_n^{-1}+\nu_m^{-1})
\;.
\]
Since $\nu_n^{-1}\to0$ as $n\to\infty$ it follows that 
 $\sigma_+\in D(\overline h_+)$.
 \hfill$\Box$

 \bigskip
 
 \noindent{\bf Proof of Theorem~\ref{c3.2.2}}$\;$
 It follows from  Lemma~\ref{lext3.137}  that 
 $D(\overline H_+)+\spann \sigma_+\subseteq D(\overline h_+)$.
Then since $h_{0,1}(\varphi)=h_+(\varphi)$ for all  $\varphi\in C_c^\infty(0,\infty)$
one deduces  that   $h_{0,1}(\Phi)=\overline h_+(\Phi)$ for all  $\Phi\in D(\overline H_+)+\spann \sigma_+$.
Hence $h_{0,1}=\overline h_+$ and $H_{0,1}=H_F$.
The boundary condition follows since $\varphi\in D(H_{0,1})$ requires $B_+(\sigma_+,\varphi)=0$.
But one has $B_+(\sigma_+,\varphi)= (c\, \varphi')(0_+)$.

The form $h_+$ is given by (\ref{e4.2.1}) with the integral restricted to the half-line. 
Then the  closure $\overline h_+$ is a Dirichlet form by standard estimates. 
Therefore  the semigroup generated by $H_{0,1}$ is submarkovian.
Now suppose that the semigroup $S^{\alpha,\beta}$ generated by $H_{\alpha,\beta}$ with $\alpha\neq0$
is submarkovian. 
Then $h_{\alpha,\beta}$ is a Dirichlet form.
In particular if $\varphi\in D(h_{\alpha,\beta})$ is positive then $R\wedge\varphi\in D(h_{\alpha,\beta})$
for all $R>0$ 
and 
\begin{equation}
h_{\alpha,\beta}(R\wedge\varphi)\leq h_{\alpha,\beta}(\varphi)
\;.
\label{edf}
\end{equation}
But it follows from the definition of $D_{\alpha,\beta}=D(H_{\alpha,\beta})$ that $\beta\,\sigma_+-\alpha\,\tau_+\in D_{\alpha,\beta}\subseteq D(h_{\alpha,\beta})$.
Moreover $\sigma_+\in D(\overline h_+)\subseteq D(h_{\alpha,\beta})$ by Lemma~\ref{lext3.137}.
Since $\alpha\neq0$ one must then  have $\tau_+\in D(h_{\alpha,\beta})$.
In particular $h_{\alpha,\beta}(\tau_+)<\infty$.
Now we can apply (\ref{edf}) with $\varphi=\tau_+$.
Since   $R\wedge \tau_+\in D(\overline h_+)\subseteq D(h_{\alpha,\beta})$ one has 
\[
\overline h_+(R\wedge \tau_+)=h_{\alpha,\beta}(R\wedge\tau_+)\leq h_{\alpha,\beta}(\tau_+)
\;.
\]
But $\tau_+=\nu_+$ and $\tau_+'=-c^{-1}$ on $\langle0,1]$.
Therefore 
\[
\int_{S_R}dx\,c(x)^{-1}= \int_{S_R}dx\,c(x)\,|\tau_+'(x)|^2\leq \overline h_+(R\wedge\tau_+)= h_{\alpha,\beta}(R\wedge\tau_+)
\leq h_{\alpha,\beta}(\tau)<\infty
\]
where $S_R=\{x\in\langle0,1]:\tau_+(x)\leq R\}$.
But $\nu_+\not\in L_\infty(0,1)$.
So  the supremum of the left hand side over $R$ is infinite.
This is  a contradiction so $h_{\alpha, \beta}$ cannot be a Dirichlet form.~\hfill$\Box$

\begin{remarkn}\label{rext2.10}
If  $\nu_+\in L_2(0,1)$ but $\nu_+\not\in L_\infty(0,1)$ then the form $h_{0,1}$ of the Friedrichs extension $H_{0,1}$
has the integral representation $h_{0,1}(\varphi)=\int^\infty_0 c\,|\varphi'|^2$ for all $\varphi\in D(h_{0,1})$
and the operator domain $D(H_{0,1})$ is distinguished by the boundary condition $(c\,\varphi')(0_+)=0$.
But the situation is quite different for the forms $h_{\alpha, \beta}$ and operators $H_{\alpha,\beta}$ if $\alpha\neq0$.
Then there is no analogue of the form representation  nor of the  boundary condition.
For example, the foregoing argument establishes that $\tau_+\in D(h_{\alpha,\beta})$.
But then
\[
\int^\infty_0 c\,|\tau_+'|^2\geq \int^1_0 c\,|\nu_+'|^2=\int^1_0 c^{-1}=\nu_+(0)=\infty
\;.
\]
Moreover, if $\Phi=\varphi+a\,\sigma_+-b\,\tau_+\in D(H_{\alpha,\beta})$ with $\varphi\in D(\overline H)$  one has
the boundary condition $a\,\alpha=b\,\beta$.
Then $b=(c\,\Phi')(0_+)$ but one cannot identify $a$ in terms of the value of $\Phi$ and its derivatives at the origin.
\end{remarkn}

Although the semigroups $S^{\alpha, \beta}$ generated by the $H_{\alpha, \beta}$ with $\alpha\neq 0$ cannot be submarkovian
we next argue that they are  positive.
First since $H_+$ has deficiency indices $(1,1)$ there is, for each $\gamma>0$,
a unique, up to a multiplicative factor, $L_2$-solution of the deficiency equation $(\gamma I+H_+^*)\eta=0$.
Therefore there is a unique positive, decreasing, normalized, $L_2$-solution $\eta_\gamma$.
More precisely, $\eta_\gamma$ is non-negative and non-increasing.
To establish this note that the equation 
 has the explicit form $(c\,\eta')'=\gamma\,\eta\in L_2(0,\infty)$ and this implies that $\eta'$ is continuous.
 Hence $\eta$ is a locally $C^1$-function.
Now  suppose $\eta(a)=0$ for some $a>0$. 
Then $(c\,\eta\,\eta')(a)=0$ and it follows from (\ref{eext2.20}), with $\varphi=\eta$ and $H^*_+\eta=-\gamma\eta$, that
\[
2^{-1}c\,(\eta^2)'(b)=(c\,\eta\,\eta')(b)=\int^b_a(\gamma\,|\eta|^2+\,c\,|\eta'|^2)
\]
for all $b\geq 0$.
In particular $(\eta^2)'(b)\geq0$ for all $b\geq a$.
But since $\eta$ is square integrable $\eta^2$ cannot be monotonically increasing.
Therefore $\eta(b)=0$ for all $b\geq a$.
If, however, $b\leq a$ then  $(\eta^2)'(b)\leq0$.
Therefore $\eta$ is either non-negative and non-increasing, or non-positive and non-decreasing.

\begin{prop}\label{p3.2.4}
Assume  $\nu_+\in L_2(0,1)$ but $\nu_+\not\in L_\infty(0,1)$. 
Then the semigroup $S^{\alpha, \beta}$ generated by the self-adjoint extension $H_{\alpha,\beta}$ of $H_+$
is positive.
\end{prop}
\proof\
If $\alpha=0$ then the semigroup is submarkovian by Theorem~\ref{c3.2.2}. 
Therefore it remains to prove positivity for $\alpha\neq0$.
Let $P_\gamma$ denote the one-dimensional orthogonal projection on $L_2(0,\infty)$ with range $\eta_\gamma$.
Since~$\eta_\gamma$ is  positive (non-negative)  the projection $P_\gamma$ is a positive operator, i.e.\ it maps positive functions into positive functions.

Next choose $\gamma$ such that $H_{\alpha, \beta}\geq \gamma I>0$.
Since $H_{0,1}=H_F\geq H_{\alpha,\beta}$ one also has $H_{0,1}+\gamma I>0$.
Then by Krein's theory of lower semibounded extensions there is a $\kappa(\gamma)\geq 0$ such that 
\begin{equation}
(\gamma I+H_{\alpha,\beta})^{-1}=(\gamma I+H_{0,1})^{-1}+\kappa(\gamma)\,P_\gamma
\label{eext3.1}
\end{equation}
(see \cite{Far}, Theorem~15.1).
But $H_{0,1}$ is a submarkovian extension so $(\gamma I+H_{0,1})^{-1}$ is a positive operator.
In addition $P_\gamma$ is a positive operator and  $\kappa(\gamma)\geq 0$.
It follows immediately that $(\gamma I+H_{\alpha,\beta})^{-1}$ is positive for all  large $\gamma$.
Then $S^{\alpha, \beta}$ is positive by the Trotter product formula.
\hfill$\Box$

\bigskip

\noindent{\bf Proof of Theorem~\ref{tnsm1.1}.\ref{tnsm1.1-2}} $\;$This follows from Theorem~\ref{c3.2.2} and Proposition~\ref{p3.2.4}.
\hfill\hfill$\Box$

\begin{remarkn}\label{rext2.11} The  representation (\ref{eext3.1})  gives information about the possible extension of the resolvents 
$(\gamma I+H_{\alpha,\beta})^{-1}$ to the $L_p$-spaces.
Since $H_{0,1}$ is a submarkovian generator it follows from (\ref{eext3.1}) that $(\gamma I+H_{\alpha,\beta})^{-1}$
extends to a bounded operator on $L_p(0,\infty)$, with $p\in\langle2,\infty\rangle$, if and only if  $P_\gamma$ extends
to a bounded operator.
Now $\eta_\gamma\in L_2(0,\infty)$ by definition.
Therefore $\eta_\gamma\in L_q(0,1)$ for all $q\in[1,2]$.
But since $c$ is bounded away from zero on $[1,\infty\rangle $ it follows that 
$\eta_\gamma\in L_r(1,\infty)$ for all $r\in[1,\infty\rangle$ by standard strong ellipticity estimates.
Therefore $\eta_\gamma\in L_q(0,\infty)$ for all $q\in[1,2]$.
Hence one concludes  that $(\gamma I+H_{\alpha,\beta})^{-1}$
extends to a bounded operator on $L_p(0,\infty)$ if and only if $\eta_\gamma\in L_p(0,1)$.

The small $x$ behaviour of $\eta_\gamma$ can, however,  be deduced from integration of the deficiency equation 
$(\gamma I+H^*)\eta_\gamma=0$.
Let $x_0\in \langle0,1]$.
After two integrations one finds
\[
\eta_\gamma(x)=\eta_\gamma(x_0)-(c\,\eta_\gamma')(x_0)\int^{x_0}_xds\,c(s)^{-1}+\gamma\int^{x_0}_xds\,c(s)^{-1}\int^{x_0}_sdt\,\eta_\gamma(t)
\]
and this leads to the estimate
\[
|\eta_\gamma(x)-\eta_\gamma(x_0)+(c\,\eta_\gamma')(x_0)\,\nu_{x_0}(x)|\leq \gamma\,\|\eta_\gamma\|_2\,|x-x_0|^{1/2}\nu_{x_0}(x)
\]
where $\nu_{x_0}(x)=\int^{x_0}_xds\,c(s)^{-1}$.
It follows  that  $\eta_\gamma\in L_p(0,1)$ if and only if $\nu_+\in L_p(0,1)$.
Therefore  $(\gamma I+H_{\alpha,\beta})^{-1}$
extends to a bounded operator on $L_p(0,\infty)$ with $p\in\langle2,\infty\rangle$ if and only if $\nu_+\in L_p(0,1)$.

The question whether the extension of the resolvent to  $L_p$ is the resolvent of the generator of an $L_p$-continuous 
semigroup seems more complicated.
\end{remarkn}

It remains to prove Theorem~\ref{tnsm1.1}.\ref{tnsm1.1-3}.
The assumption  $\nu_+\in L_\infty(0,1)$  corresponds to integrability  of $c^{-1}$ at the origin.
Therefore  it is natural to reparametrize and replace $\nu_+$ by $\hat\nu$ where   
\[
\hat\nu(x)=\nu_+(0)-\nu_+(x)=\int^x_0c^{-1}
\]
and then to  replace $\tau_+=\nu_+\,\sigma_+$ by $\hat \tau=\hat \nu\,\sigma_+$.
Since $\hat \tau=\nu_+(0)\,\sigma_+-\tau_+$ these replacements make no essential difference to the characterization of $D(H_+^*)$ given by
Proposition~\ref{p3.2.1}.
Now one has 
$D(H_+^*)=D(\overline H_+)+\spann \sigma_+ +\spann \hat \tau$.
The advantage of the reparametrization is that if $\Phi=\varphi+a\,\sigma_++b\,\hat \tau\in  D(H_+^*)$ with $\varphi\in D(\overline H_+)$ 
then  $a=\Phi(0_+)$ and $b=(c\,\Phi')(0_+)$.
This follows since $\varphi(0_+)=0=(c\,\varphi')(0_+)$, by Proposition~\ref{pext2.1}.\ref{pext2.1-2}, and $\sigma(0_+)=1$, $(c\,\sigma_+')(0_+)=0$, $\hat \tau(0)=0$ and $(c\,\hat \tau')(0_+)=1$
by definition.
Now we modify accordingly the definition of the self-adjoint extensions of $H_+$.

Let $(\alpha,\beta)\in\Ri^2\backslash(0,0)$.
Define 
\[
\widehat D_{\alpha,\beta}=D(\overline H_+)+\spann(\beta \,\sigma_++\alpha\,\hat \tau)
\;.
\]
Then define the self-adjoint extension $\widehat H_{\alpha,\beta}$ of $H$ as the restriction of $H_+^*$ to $\widehat D_{\alpha,\beta}$. 
Again one has  $\widehat H_{\alpha,\beta}=\widehat H_{\alpha',\beta'}$ if and only if $\alpha\,\beta'=\beta\,\alpha'$.
Further let $\hat h_{\alpha,\beta}$ denote the quadratic form corresponding to  $\widehat H_{\alpha,\beta}$.

One can again compute the Friedrichs  extension of $H_+$.

\begin{prop}\label{c3.3.1}
Assume  $\nu_+\in L_\infty(0,1)$. 
Then  $\overline h_+=\hat h_{1,0}$, the Friedrichs  extension $H_F$ of $H_+$ is equal to $\widehat H_{1,0} $
 and 
\[
D(H_F)=\{\varphi\in D(H_+^*):  \varphi(0_+)=0\}
\;.
\]
\end{prop}
\proof\ First one has  $\hat h_{1,0}\supseteq \overline h_+$ and consequently $D(\hat h_{1,0})\supseteq D(\overline h_+)$.
Secondly,  $\hat\tau\in D(\overline h_+)$.
This follows by standard approximation techniques since $\hat \tau$ has compact support, is absolutely continuous and $\hat \tau(0_+)=0$.
But then  $D(\widehat H_{1,0})=D(\overline H_+)+\spann\hat\tau\subseteq D(\overline h_+)$.
Since $D(\widehat H_{1,0})$  is a core of $\hat h_{1,0}$ it  follows that $D(\overline h_+)$ is also a core of $\hat h_{1,0}$.
Therefore $\hat h_{1,0}=\overline h_+$ and $\widehat H_{1,0}=H_F$.

The boundary condition follows since if $\Phi=\varphi+a\,\sigma_++b\,\hat\tau\in  D(H_+^*)$ with $\varphi\in D(\overline H_+)$ 
then $\Phi\in D(\widehat H_{1,0})$ if and only if $a=0$.
But it follows from the  reparametrization chosen above that    $a= \Phi(0_+)$.
\hfill$\Box$

\bigskip

The analysis of the remaining self-adjoint extensions of $H_+$ is  in terms of the corresponding quadratic forms.
First we observe that there is a unique form domain.

\begin{prop}\label{pnsm3.11} If $\beta\neq 0$ then 
$D(\hat h_{\alpha,\beta})=D(\overline h_+)+\spann\sigma_+=D(\hat h_{0,1}).$
\end{prop}
\proof\
First, since $\hat h_{\alpha, \beta}\supseteq \overline h_+$ one has $D(\overline h_+)\subseteq D(\hat h_{\alpha,\beta})$.
But $\hat\tau\in D(\overline h_+)$, as observed in the proof of Proposition~\ref{c3.3.1}.
Moreover,  since $D(\widehat H_{\alpha,\beta})\subseteq D(\hat h_{\alpha, \beta})$ it follows that $a\,\sigma_++b\,\hat \tau\in D(\hat h_{\alpha,\beta})$
for all $a, b$ satisfying $a\,\alpha=b\,\beta$.
Since $\beta\neq0$ one deduces that $\sigma_+\in D(\hat h_{\alpha,\beta})$.
Therefore  $D(\hat h_{\alpha,\beta})\supseteq D(\overline h_+)+\spann\sigma_+$.

Secondly, we establish the converse inclusion.
The proof begins by observing that 
$D(\widehat H_{\alpha,\beta})$ is a core of $\hat h_{\alpha,\beta}$.
Thus if $\Phi\in D(\hat h_{\alpha,\beta})$ there is a sequence $\Phi_n\in D(\widehat H_{\alpha,\beta})$ which converges to $\Phi$ in the $D(\hat h_{\alpha,\beta})$-graph norm.
But $\Phi_n=\varphi_n+a_n\,\sigma_+ +b_n\,\hat \tau$ with $\varphi_n\in D(\overline H_+)$ and $a_n, b_n\in \Ri$ satisfying $a_n\,\alpha=b_n\,\beta$.
Moreover,   $a_n=\Phi_n(0_+)$  and $b_n=(c\,\Phi_n')(0_+)$ by the new choice of parameters. 
Therefore
\[
|a_n-a_m|=|\Phi_n(0_+)-\Phi_m(0_+)|=|B_+(\Phi_n-\Phi_m,\hat\tau)|=|(H_+^*(\Phi_n-\Phi_m),\hat\tau)-((\Phi_n-\Phi_m),H_+^*\hat\tau)|
\]
where the second equality follows because $\hat\tau(0_+)=0$ and $(c\,\hat\tau')(0_+)=1$.
But $( H^*\hat\tau)(x)$ is bounded with support in the interval $[1,2]$.
Hence $|((\Phi_n-\Phi_m),H^*\hat \tau)|\leq \kappa\,\|\Phi_n-\Phi_m\|_2$ for some $\kappa>0$.
Moreover, $\hat \tau\in D(\overline h_+)\subseteq D(\hat h_{\alpha,\beta})$ by the proof of Proposition~\ref{c3.3.1}.
Hence 
\[
|(H_+^*(\Phi_n-\Phi_m),\hat \tau)|=|(\widehat H_{\alpha,\beta}(\Phi_n-\Phi_m),\hat \tau)|\leq \hat h_{\alpha,\beta}(\Phi_n-\Phi_m)^{1/2}\,\hat h_{\alpha,\beta}(\hat \tau)^{1/2}
\;.
\]
Therefore
\[
|a_n-a_m|\leq \hat h_{\alpha,\beta}(\Phi_n-\Phi_m)^{1/2}\,\hat h_{\alpha,\beta}(\tau)^{1/2}+\kappa\,\|\Phi_n-\Phi_m\|_2
\]
and   the $a_n$ must converge to a limit $a$.
Similarly, 
\[
|b_n-b_m|=|(c\,\Phi'_n)(0_+)-(c\,\Phi'_m)(0_+)|=|B_+(\Phi_n-\Phi_m,\sigma_+)|
\]
 because $\sigma_+(0_+)=1$ and $(c\,\sigma_+')(0_+)=0$.
 But then\[
|b_n-b_m|\leq \hat h_{\alpha,\beta}(\Phi_n-\Phi_m)^{1/2}\,\hat h_{\alpha,\beta}(\sigma_+)^{1/2}+\kappa\,\|\Phi_n-\Phi_m\|_2
\]
and   the $b_n$ must converge to a limit $b$.
One automatically has $a\,\alpha=b\,\beta$.

Next it follows that $\varphi_n=\Phi_n-a_n\,\sigma_+-b_n\hat\tau$ is $L_2$-convergent to a limit $\varphi$.
But 
\begin{eqnarray*}
\overline h_+(\varphi_n-\varphi_m)&=&\hat h_{\alpha,\beta}(\varphi_n-\varphi_m)\\[5pt]
&\leq& 2\,\hat h_{\alpha,\beta}(\Phi_n-\Phi_m)+4\,|a_n-a_m|\,\hat h_{\alpha,\beta}(\sigma_+)
+4\,|b_n-b_m|\,\hat h_{\alpha,\beta}(\hat\tau)
\end{eqnarray*}
so the $\varphi_n$ are convergent in the $D(\overline h_+)$-graph norm and one has $\varphi\in D(\overline h_+)$.
Therefore $\Phi=\varphi+a\,\sigma_++b\,\hat\tau\in D(\overline h_+)+\spann\sigma_+$.
Hence  $D(\hat h_{\alpha,\beta})\subseteq D(\overline h_+)+\spann\sigma_+$.
\hfill$\Box$

\bigskip

Finally one can express the forms $\hat h_{\alpha,\beta}$ in terms of $\hat h_{0,1}$ on their
common domain in a classical manner.

\begin{prop}\label{tsm3.3}
If $\nu_+\in L_\infty(0,1)$ then 
$D(\hat h_{\alpha,\beta})=D(\hat h_{0,1})$ for all $\beta\in\Ri\backslash \{0\}$
and
\begin{equation}
\hat h_{\alpha,\beta}(\varphi)=\hat h_{0,1}(\varphi)+\alpha\,\beta^{-1}\,|\varphi(0_+)|^2
\label{fid1}
\end{equation}
for all $\varphi\in D(h_{0,1})$.

The self-adjoint extension   $\widehat H_{\alpha,\beta}$ is the restriction of $H_+^*$ to the domain
\[
D(\widehat H_{\alpha,\beta})=\{\varphi\in D(H_+^*): \beta\, (c\,\varphi')(0_+)=\alpha\,\varphi(0_+)\}\;.
\]
The operators   $\widehat H_{\alpha,\beta}$ generate positive semigroups on $L_2(0,\infty)$
which are submarkovian if and only if 
 $\alpha\,\beta^{-1}\geq 0$.
\end{prop}
\proof\
The identity of the domains is established in Proposition~\ref{pnsm3.11}.
Next we establish the relation between the forms.

First suppose $\Phi=\varphi+a\,\sigma_++b\,\hat\tau$ with $\varphi\in D(\overline H_+)$.
Then  $\Phi_1=\varphi+a\,\sigma_+\in D(\widehat H_{0,1})\subseteq D(\hat h_{0,1})$.
Moreover, $\hat\tau\in D(\overline h_+)$ as observed in the proof of Proposition~\ref{c3.3.1}. 
Therefore $\Phi\in  D(\hat h_{0,1})$.
Then one calculates that 
\begin{eqnarray*}
(\Phi,H_+^*\Phi)&=&\hat h_{0,1}(\Phi_1)+b\,(\hat\tau,H_+^*\Phi_1)+b\,(\Phi_1,H_+^*\hat\tau)+b^2\,(\hat\tau,H_+^*\hat\tau)\\[5pt]
&=&
\hat h_{0,1}(\Phi_1)+
b\,(\hat\tau,H_+^*\Phi_1)+b\,(H_+^*\Phi_1,\hat\tau)+b^2\,(\hat\tau,H_+^*\hat\tau)+b\,B_+(\Phi_1,\hat\tau)
\;.
\end{eqnarray*}
But 
\[
(\hat\tau,H_+^*\Phi_1)=(\hat\tau,\widehat H_{0,1}\Phi_1)=\hat h_{0,1}(\hat\tau,\Phi_1)
\]
since $\hat\tau\in  D(\hat h_{0,1})$.
Similarly  $(H_+^*\Phi_1,\hat\tau)=\hat h_{0,1}(\Phi_1,\hat\tau)$ and 
$(\hat\tau,H_+^*\hat\tau)=\hat h_{0,1}(\hat\tau)$.
Combining these identities gives 
\begin{eqnarray*}
(\Phi,H_+^*\Phi)&=&\hat h_{0,1}(\Phi_1)+b\,\hat h_{0,1}(\hat\tau,\Phi_1)
+b\,\hat h_{0,1}(\Phi_1,\hat\tau)+b^2\,\hat h_{0,1}(\hat\tau)+b\,B_+(\Phi_1,\hat\tau)\\[5pt]
&=&\hat h_{0,1}(\Phi)+b\,B_+(\Phi_1,\hat\tau)=\hat h_{0,1}(\Phi)+a\,b
\end{eqnarray*}
where the last step uses the identification $B_+(\Phi_1,\hat\tau)=a\,B(\sigma_+,\hat\tau)=a$.
If, however, one places the restriction $a\,\alpha=b\,\beta$ on $a$ and $b$  then $\Phi\in D(\widehat H_{\alpha,\beta})$ and $H_+^*\Phi=\widehat H_{\alpha,\beta}\Phi$.
Therefore 
\[
\hat h_{\alpha,\beta}(\Phi)=\hat h_{0,1}(\Phi)+a\,b=\hat h_{0,1}(\Phi)+\alpha\,\beta^{-1}\,a^2
\]
for all $\Phi\in D(\widehat H_{\alpha,\beta})$ and then by closure for all  $\Phi\in D(\hat h_{\alpha,\beta})$.
But  $a=\Phi(0_+)$ so this etablishes the relation (\ref{fid1}).

The boundary condition for $\Phi\in D(\widehat H_{\alpha,\beta})$ is given by $a\,\alpha=b\,\beta$
but the parametrization was chosen such that $a=\Phi(0_+)$ and $b=(c\,\Phi')(0_+)$.

Finally if $\Phi=\varphi+a\,\sigma_+$ with $\varphi\in D(\overline H_+)$ then since  $(c\,\Phi')(0_+)=0$ one computes that 
\[
\hat h_{0,1}(\Phi)=(\Phi,H_+^*\Phi)=-(\Phi,(c\,\Phi')')=\int^\infty_0dx\, c(x)\,|\Phi'(x)|^2\;.
\]
Therefore 
\[
\hat h_{\alpha,\beta}(\Phi)=\int^\infty_0dx\, c(x)\,|\Phi'(x)|^2+(\alpha\,\beta^{-1})|\Phi(0_+)|^2
\]
and  the positivity and submarkovian properties follow immediately by application of the
well known Beurling--Deny criteria (see, for example, \cite{RS4}, pages~209--212).
\hfill$\Box$

\bigskip

\noindent{\bf Proof of Theorem~\ref{tnsm1.1}.\ref{tnsm1.1-3}}$\;$ 
This follows directly from Proposition~\ref{c3.3.1}  and Theorem~\ref{tsm3.3}.
These latter results give a complete description of the self-adjoint extensions
of the operator $H$  for  $\nu_+\in L_\infty(0,1)$.
The Friedrichs extension $\widehat H_{1,0}$ corresponds to Dirichlet boundary conditions $\varphi(0_+)=0$,  the extension
$\widehat H_{0,1}$ to Neumann boundary conditions $(c\,\varphi')(0_+)=0$ and the other extensions to Robin boundary conditions.
The latter are positive-definite if $\alpha\,\beta^{-1}\geq 0$.
\hfill$\Box$

\begin{remarkn} \label{rext2.12} 
It is straightforward to establish that if    $\alpha\,\beta^{-1}<0$ then the extensions  $H_{\alpha,\beta}$ of $H$ have a simple negative
eigenvalue.
Therefore these extensions are no longer contractive on $L_2(0,\infty)$ and certainly not submarkovian.
\end{remarkn}

\begin{remarkn}\label{rext2.13} 
Theorem~\ref{tnsm1.1} identifies  two distinct  cases in which 
 $H_+$ has a unique submarkovian extension, the Friedrichs extension $ H_{+F}$.
In both cases $\nu_+\not\in L_\infty(0,1)$.
 Moreover, in both cases
  $D( H_{+F})=D(\overline H_+)+\spann\sigma_+$ and 
$H_{+F}(\varphi+a\,\sigma_+)=\overline H_+\varphi-a\,(c\,\sigma_+')'$
for all $\varphi\in D(\overline H_+)$ and $a\in\Ri$.
The distinction between the cases occurs because $\nu_+\not\in L_2(0,1)$  implies that 
$\sigma_+\in D(\overline H_+)$.
Therefore  $ H_{+F}=\overline H_+$ and $H_+$ is essentially self-adjoint.
But if $\nu_+\in L_2(0,1)$ then $\sigma_+\not\in D(\overline H_+)$ and $ H_{+F}$ is a strict extension of $\overline H_+$.

\end{remarkn}

\section{The line}\label{S3}

In this section we prove Theorem~\ref{text4.1} by applying the results of Section~\ref{S2}
to $H_+=H|_{C_c^\infty(0,\infty)}$ on $L_2(0,\infty)$,  $H_-=H|_{C_c^\infty(-\infty,0)}$ 
on $L_2(-\infty,0)$ and $H_0=H_-\oplus H_+$ on $L_2(\Ri)$.

Let $\sigma\in C_c^\infty(\Ri)$ satisfy $0\leq \sigma\leq1$, $\supp\sigma\subseteq [-2,2]$ and  $\sigma=1$ on $\langle-1,1\rangle$.
Define $\sigma_+$ by $\sigma_+(x)=0$ if $x<0$ and $\sigma_+(x)=\sigma(x)$ if $x\geq 0$
and set $\sigma_-=\sigma-\sigma_+$.
Further define $\tau_\pm$ by $\tau_\pm(x)=\nu_\pm(\pm \,x)\,\sigma_\pm(x)$.
Thus $\tau_+$ has support in $[0,2]$ and $\tau_-$ has support in  $[-2,0]$.
Now one can characterize the self-adjoint extensions of $H_+$ with the aid of the functions $\sigma_+, \tau_+$ exactly as in Section~\ref{S2}
and the extensions of $H_-$ with the aid of $\sigma_-,\tau_-$ in an analogous fashion.

\smallskip

\noindent{\bf Proof of Theorem~\ref{text4.1}.\ref{text4.1-1}}$\;$
There are two   cases to be considered: 1. $\nu_+\wedge \nu_-\not\in L_2(0,1)$,  and 2. $\nu_+\wedge \nu_-\in L_2(0,1)$
and $\nu_+\vee \nu_-\not\in L_2(0,1)$.

Assume $\nu_+\wedge \nu_-\not\in L_2(0,1)$.
Then $H_+$ is essentially self-adjoint on $L_2(0,\infty)$  and $H_-$ is essentially self-adjoint on $L_2(-\infty,0)$ by Theorem~\ref{tnsm1.1}.\ref{tnsm1.1-1}.
Hence $\overline H_0=\overline H_-\oplus\overline H_+$ is self-adjoint on $L_2(\Ri)$.
Since a self-adjoint operator cannot have a proper closed symmetric extension it follows that $\overline H=\overline H_0=\overline H_-\oplus\overline H_+$ is self-adjoint.
Therefore $\overline H$ must coincide with the Friedrichs extension $H_F$ of $H$ and it automatically generates a submarkovian semigroup.
Clearly this semigroup must leave $L_2(-\infty,0)$ and $L_2(0,\infty)$ invariant.

Secondly, assume $\nu_+\not\in L_2(0,1)$  and $\nu_-\in L_2(0,1)$.
Then the argument is slightly different although the conclusion is the same.
Again  $H_+$ is essentially self-adjoint on $L_2(0,\infty)$.
But then  $\sigma_+\in  D(H_+^*)=D(\overline H_+)$ 
or, as a relation on $L_2(\Ri)$, 
 $\sigma_+\in  D(\overline H)$.
Since  $\sigma\in D(H)$ it follows that $\sigma_- \in D(\overline H)$.
Now one can define a self-adjoint  extension $\widetilde  H_-$ of $H_-$  by $D(\widetilde H_-)=D(\overline H_-)+\spann \sigma_-\subseteq D(\overline H)$ and 
$\widetilde H_-(\varphi+\beta\,\sigma_-)=\overline H_-\varphi-\beta\,(c\,\sigma_-')'$ for $\varphi\in D(\overline H_-)$
and $\beta\in \Ri$.
Then $\overline H\supseteq \widetilde H_- \oplus\overline H_+ $ and since  a self-adjoint operator cannot have a proper closed symmetric extension it follows that 
$\overline H=\widetilde H_- \oplus\overline H_+$ is self-adjoint.
Then $\overline H$ must coincide with the Friedrichs extension $H_F$ of $H$ and the corresponding semigroup is submarkovian.
Clearly the semigroup leaves $L_2(-\infty,0)$ and $L_2(0,\infty)$ invariant.

Thirdly, the argument for  $\nu_+\in L_2(0,1)$  and $\nu_-\not\in L_2(0,1)$ is similar.\hfill$\Box$

\smallskip

\noindent{\bf Proof of Theorem~\ref{text4.1}.\ref{text4.1-2}}$\;$
Now we assume $\nu_+\vee \nu_-\in L_2(0,1)$.
Since $H_0\subseteq H$ and $H_0$  and $H$ are both symmetric one has 
\begin{equation}
H_0\subseteq H\subseteq H^*\subseteq H^*_0
\;.
\label{enext1}
\end{equation}
But $H_0$ is the direct sum of $H_\pm$  and as both these operators have deficiency indices $(1,1)$, by Proposition~~\ref{pext2.1}.\ref{pext2.1-2},
 the operator $H_0$ must have deficiency indices $(2,2)$.
Thus $D(\overline H_0)$ has codimension 4 in $D(H_0^*)$.
Moreover, $H_0^*=H_-^*\oplus H_+^*$ and one can compute the adjoint of $H_0$ in terms of the adjoints of the operators
$H_\pm$ on the half-lines.
But this allows one to compute the domain of $H^*$.

\begin{prop}\label{pext4.1}
If  $\nu_+\vee \nu_-\in L_2(0,1)$
then 
\begin{eqnarray}
D(H^*)&=&D(\overline H_0 )+\spann\sigma_++\spann\sigma_-+\spann(\tau_+-\tau_-)\nonumber\\[5pt]
&=&D(\overline{H})+\spann(\sigma_+-\sigma_-)+\spann(\tau_+-\tau_-)
\;.
\label{enext2.0}
\end{eqnarray}
\end{prop}
\proof\
It follows from Proposition~\ref{p3.2.1} applied to $H_\pm$  that 
\begin{equation}
D(H_0^*)=D(\overline H_0 )+\spann\sigma_++\spann\sigma_-+\spann\tau_++\spann\tau_-
\;.
\label{enext2}
\end{equation}
Now introduce the boundary form $\varphi,\psi\in D(H_0^*)\mapsto B_0(\varphi,\psi)$ by
\[
B_0(\varphi,\psi)=(H_0^*\varphi,\psi)-(\varphi,H_0^*\psi)
\;.
\]
Then it follows from (\ref{enext1})  that 
\begin{equation}
D(H^*)=\{\varphi\in D(H_0^*):B_0(\varphi,\psi)=0 \mbox{ for all } \psi\in D(H)\}
\label{enext3}
\end{equation}
(see \cite{Far}, lemma on page 86).
Now one can compute $D(H^*)$ by use of (\ref{enext2}).

First  $B_0(\varphi,\psi)=0$ for all $\varphi\in D(\overline H_0 )$ and $\psi\in D(H)$.
Secondly
\[
B_0(\sigma_+,\psi)=0=B_0(\sigma_-,\psi)
\]
for all $\psi\in C_c^\infty(\Ri)$ by direct calculation.
Thirdly,
\[
B_0(\tau_+,\psi)=(c\,\nu_+'\,\psi)(0)-(\nu_+\,c\,\psi')(0)=-\psi(0)=B_0(\tau_-,\psi)
\]
since $(\nu_+\,c\,\psi')(0)=0=(\nu_-\,c\,\psi')(0)$ for all $\psi\in C_c^\infty(\Ri)$.
The latter relations follow because   the assumption $c(0)=0$  implies that $(c\,\psi')(0)=0$.
Hence
$(c\,\psi')(x)=\int^x_0 ds\,(c\,\psi')'(s)$  and this gives
 an estimate $|(c\,\psi')(x)|\leq |x|^{1/2}\,\|H\psi\|_2$.
As the $\nu_\pm$ are square integrable near the origin one has $|x|^{1/2}\nu_\pm(x)\to0$ as $x\to 0$ by Lemma~\ref{lext2.11}.
 Therefore $B_0(\tau_+-\tau_-,\psi)=0$ for all $\psi\in D(H)$. 
 But there is a $\psi\in D(H)$ such that 
 $B_0(\tau_++\tau_-,\psi)\neq0$.
Consequently one concludes from (\ref{enext3}) that
 \[
 D(H^*)= D(\overline H_0)+\spann\sigma_++\spann\sigma_-+\spann(\tau_+-\tau_-)
 \]
 which is the first statement of the proposition.
 But $\sigma_++\sigma_-=\sigma\in C_c^\infty(\Ri)\subseteq D(\overline H)$.
 Hence $D(\overline H_0)+\spann(\sigma_++\sigma_-)\subseteq D(\overline H)$.
Therefore 
 \begin{eqnarray*}
 D(H^*)&=& D(\overline H_0)+\spann(\sigma_++\sigma_-)+\spann(\sigma_+-\sigma_-)+\spann(\tau_+-\tau_-)\\[5pt]
&\subseteq& D(\overline H)+\spann(\sigma_+-\sigma_-)+\spann(\tau_+-\tau_-)\subseteq D(H^*)
 \end{eqnarray*}
which gives the second statement of the proposition.
\hfill$\Box$

\bigskip

Note that $B_0(\tau_+,\sigma_+)=1=B_0(\tau_-,\sigma_-)$.
Hence, under the assumptions of the proposition, one cannot have $\sigma_+-\sigma_-,\tau_+-\tau_-\in D(\overline H)$.
Therefore $D(\overline H)$ has codimension 2 in $D(H^*)$, i.e.\ $H$ has deficiency indices $(1,1)$.
Moreover, since $\sigma\in D(H)$ one has 
\[
D(H^*)=D(\overline{H})+\spann\sigma_++\spann(\tau_+-\tau_-)=D(\overline{H})+\spann\sigma_-+\spann(\tau_+-\tau_-)
\;.
\]

The self-adjoint extensions  $H_{\alpha,\beta}$  of $H$ are given  for  $(\alpha,\beta)\in \Ri^2\backslash(0,0)$
by 
\[
D(H_{\alpha,\beta})=D(\overline H)+\spann\Big(\beta\,(\sigma_+-\sigma_-)-\alpha\,(\tau_+-\tau_-)\Big)
\]
and $H_{\alpha,\beta}\Phi=H^*\Phi$ for $\Phi\in D(H_{\alpha,\beta})$.
This definition is the direct analogue of the definition on the half-line given in Proposition~\ref{p3.2.2}.
Again  $H_{\alpha,\beta}=H_{\alpha',\beta'}$ for all pairs with $\alpha\,\beta'=\alpha'\,\beta$.
In terms of the boundary form $\varphi,\psi\in D(H^*)\mapsto B(\varphi,\psi)$ associated with $H$ one has
\[
D(H_{\alpha,\beta})=\{\varphi\in D(H^*): B(\beta\,(\sigma_+-\sigma_-)-\alpha\,(\tau_+-\tau_-),\varphi)=0\}
\;.
\]
Now consider the extension with $\alpha=0$ and $\beta=1$.
Then $D(H_{0,1})=D(\overline H)+\spann(\sigma_+-\sigma_-)$.
But if $\nu_+\not\in L_\infty(0,1)$ then $\sigma_+\in D(\overline{h})$ by Lemma~\ref{lext3.137}
applied to $h$ instead of $h_+$. 
Since $\sigma=\sigma_++\sigma_-\in C^\infty_c(\Ri)$ it follows that $\sigma_-\in D(\overline{h})$. 
Thus $D(H_{0,1})\subseteq D(\overline h)$. 
Therefore $h_{0,1}=\overline h$ and 
 $H_{0,1}$ is the Friedrichs extension $H_F$ of $H$.
A similar conclusion is valid if $\nu_-\not\in L_\infty(0,1)$.

Therefore if  $\nu_+\wedge\nu_-\not\in L_\infty(0,1)$  the operator $H_{0,1}$ generates a submarkovian semigroup
$S$ which  leaves the subspaces $L_2(-\infty,0)$ and $L_2(0,\infty)$ invariant.
This establishes the first part of Theorem~\ref{text4.1}.\ref{text4.1-2}.

Note that if $\Phi=\varphi+a\,(\sigma_+-\sigma_-)-b\,(\tau_+-\tau_-)$ with $\varphi\in D(\overline H)$
then 
\[
(c\,\Phi')(0_+)=\lim_{x\to0_+}(c\,\Phi')(x)=-b\lim_{x\to0_+}(c\,\tau_+')(x)=b
\]
 and 
\[
(c\,\Phi')(0_-)=\lim_{x\to0_-}(c\,\Phi')(x)=b\lim_{x\to0_-}(c\,\tau_-')(x)=-b
\;.
\]
Therefore 
\[
b=\Big((c\,\Phi')(0_+)-(c\,\Phi')(0_-)\Big)/2
\;.
\]
In particular if $\Phi\in D(H_{0,1})$ then $b=0$ and  the extension is characterized by the boundary condition 
$(c\,\Phi')(0_+)=(c\,\Phi')(0_-)$ which links the left and right half-lines.

It remains to  prove that under the assumption of the second statement of Theorem~\ref{text4.1} there are no other
submarkovian extensions.
The key step in the proof is the  identification of the corresponding form domains.

\begin{prop}\label{pext3.11}
Assume $\nu_+\vee\nu_-\in L_2(0,1)$ but $\nu_+\vee\nu_-\not\in L_\infty(0,1)$.
If $\alpha\neq 0$ then 
$D(h_{\alpha,\beta})=D(\overline h)+\spann(\tau_+-\tau_-)=D(h_{1,0}).$
\end{prop}
\proof\
The proof is very similar to the proof of Proposition~\ref{pnsm3.11}.

First, since $ h_{\alpha, \beta}\supseteq \overline h$ one has $D(\overline h)\subseteq D(h_{\alpha,\beta})$.
But $\sigma_\pm \in D(\overline h)$, again  by Lemma~\ref{lext3.137} and the observation that 
$\sigma=\sigma_++\sigma_-\in C^\infty_c(\Ri)$.
Moreover,  since $D(H_{\alpha,\beta})\subseteq D( h_{\alpha, \beta})$ it follows that $a\,(\sigma_+-\sigma_-)-b\, (\tau_+-\tau_-)\in D( h_{\alpha,\beta})$
for all $a, b$ satisfying $a\,\alpha=b\,\beta$.
Since $\alpha\neq0$ one deduces that $\tau_+-\tau_-\in D( h_{\alpha,\beta})$.
Therefore  $ D(\overline h)+\spann(\tau_+-\tau_-)\subseteq D( h_{\alpha,\beta}) $.

Secondly,  the converse inclusion is established by a slight modification of the second part of the proof of Proposition~\ref{pnsm3.11}.
It is again dependent on the observation that $D( H_{\alpha,\beta})$ is a core of $ h_{\alpha,\beta}$.
We omit the details.\hfill$\Box$

\begin{cor}\label{cnext1}
Assume $\nu_+\vee\nu_-\in L_2(0,1)$ but $\nu_+\vee\nu_-\not\in L_\infty(0,1)$.
  If $\alpha\neq 0$ then the semigroup $S^{\alpha,\beta}$ generated $H_{\alpha,\beta}$ is neither positive nor $L_\infty$-contractive.
  \end{cor}
  \proof\
  It is necessary for positivity of $S^{\alpha,\beta}$ that $\varphi\in D(h_{\alpha,\beta})$ implies $|\varphi|\in D(h_{\alpha,\beta})$.
  This is a consequence of the first Beurling--Deny criterion (see, for example, \cite{RS4}, page~209).
  But  $\tau_+-\tau_-\in D(h_{1,0})=D(h_{\alpha,\beta})$ and by definition the $\tau_\pm$ are positive with disjoint supports.
  Therefore  $|\tau_+-\tau_-|=\tau_++\tau_-$.
  Since $\tau_\pm\not\in D(h_{1,0})$ the semigroup is not positive.
  
  The failure of $L_\infty$-contractivity is established by the argument used for the half-line (see the proof of Theorem~\ref{c3.2.2}).
  \hfill$\Box$

\bigskip

This completes the proof of Theorem~\ref{text4.1}.\ref{text4.1-2}. \hfill$\Box$

\bigskip

\noindent{\bf Proof of Theorem~\ref{text4.1}.\ref{text4.1-3}}
Assume  $\nu_+\vee\nu_-\in L_2(0,1)$.
 Then define $\nu$ by
\[
\nu(x)=\int^x_0ds\,c(s)^{-1}
\]
and $\tau=\nu\,\sigma$.
It follows readily that $\tau$ is related to the previous functions $\tau_\pm$ by a relation
\[
\tau= \gamma\,\sigma+\delta\,(\sigma_+-\sigma_-)-(\tau_+-\tau_-)
\]
with $\gamma,\delta\in\Ri$ and $\delta\geq0$.
Therefore the self-adjoint extensions $H_{\alpha,\beta}$ of $H$ can now be defined as the restrictions of $H^*$
to the domains
$D(H_{\alpha,\beta})=D(\overline H)+\spann(\beta(\sigma_+-\sigma_-)+\alpha\,\tau)$
with a typical element $\Phi\in D(H_{\alpha,\beta})$  given by
$\Phi=\varphi+a\,(\sigma_+-\sigma_-)+b\,\tau$
where $\varphi\in D(\overline H)$ and $a\,\alpha=b\,\beta$.
Therefore $\Phi(0_\pm)=\varphi(0)\pm a$, $(c\,\Phi')(0_\pm)=\pm b$ and one has $a=(\Phi(0_+)-\Phi(0_-))/2$
and $b=((c\,\Phi')(0_+)-(c\,\Phi')(0_-))/2$.
Thus $\Phi$ satisfies the boundary condition
\[
\beta\,\Big((c\,\Phi')(0_+)-(c\,\Phi')(0_-)\Big)=\alpha\,\Big(\Phi(0_+)-\Phi(0_-)\Big)
\;.
\]

A slight variation of  the previous arguments gives $h_F=h_{1,0}$ and 
$D(H_F)=D(\overline H)+\spann\tau$. 
The corresponding boundary condition is $(c\,\Phi')(0_+)=(c\,\Phi')(0_-)$.
Although the Friedrichs extension is automatically  submarkovian the corresponding
semigroup no longer leaves the subspaces $L_2(0,\infty)$ and $L_2(-\infty,0)$ invariant.
The boundary condition now links the two sides of the line.

If $\beta\neq0$ then arguing as in the proofs of  Propositions~\ref{pnsm3.11} and \ref{pext3.11} one finds
\[
D(h_{\alpha,\beta})=D(\overline h)+\spann (\sigma_+-\sigma_-)=D(h_{0,1})
\;.
\]
Then by an argument analogous to the proof of Proposition~\ref{tsm3.3} one deduces that 
\[
h_{\alpha,\beta}(\varphi)=h_{0,1}(\varphi)+\alpha\,\beta^{-1}\,|(\varphi(0_+)-\varphi(0_-)|^2/4
\;.
\]
Since $D(H_{0,1})$ corresponds to the boundary condition $(c\,\varphi')(0_+)=(c\,\varphi')(0_-)$
it follows that 
\[
h_{0,1}(\varphi)=\int_\Ri dx\,c(x)\,|\varphi'(x)|^2
\]
for all $\varphi\in D(H_{0,1})$.
Then, by closure, $h_{0,1} $ is a Dirichlet form and $H_{0,1}$ is submarkovian.

  Since $||\varphi(x)|-|\varphi(y)||\leq |\varphi(x)-\varphi(y)|$ it  follows that 
  if $\alpha\,\beta^{-1}\geq0$ then  $h_{\alpha,\beta}$  satisfies the first Beurling--Deny criterion.
  A similar argument shows that under the same restriction on $\alpha$ and $\beta$ it satisfies the second
  criterion. 
  Therefore one concludes that if $\alpha\,\beta^{-1}\geq0$ then $H_{\alpha,\beta}$ is submarkovian.
  
  Finally we note that if  $\alpha\,\beta^{-1}<0$ one can establish that $H_{\alpha,\beta}$ has a simple negative eigenvalue.
So $S^{\alpha,\beta}$ is not  contractive  on $L_2(\Ri)$ and  therefore not submarkovian.
 \hfill$\Box$

\begin{remarkn}\label{rext3.11}
If  $H$ has a unique  submarkovian extension then it is equal to the Friedrichs extension $H_F$ and 
 is given by $D(H_F)=D(\overline H)+\spann(\sigma_+-\sigma_-)$ and 
\[
H_F(\varphi+a\,(\sigma_+-\sigma_-))=\overline H\varphi-a\,((c\,\sigma_+')'-(c\,\sigma_-')')
\]
for all $\varphi\in D(\overline H)$ and $a\in\Ri$.
There are two distinct cases corresponding to the first two cases of Theorem~\ref{text4.1}.
In the first case,  $\nu_+\vee\nu_-\not\in L_2(0,1)$, one has $\sigma_\pm\in D(\overline H)$.
Therefore  $ H_F=\overline H$ and $H$ is essentially self-adjoint.
In the second case  $\sigma_\pm\not\in D(\overline H)$.
\end{remarkn}

\section{$L_1$-estimates}\label{S5}

The principal aim of this section is the proof of  Theorem~\ref{cext1.0}.
This requires a number of preliminary $L_1$-estimates which are valid
under the weaker hypothesis $c\in W^{1,1}_{\rm loc}(\Ri)$.
This is sufficient to ensure that $HC_c^\infty(\Ri)\subseteq L_1(\Ri)$.
We again begin by analyzing $H$ on the half-line.

Let $H_+=H|_{C_c^\infty(0,\infty)}$.  
Then $H_+C_c^\infty(0,\infty)\subseteq L_1(0,\infty)$
and we may  consider $H_+$ as an operator on $L_1(0,\infty)$ with domain $D(H_+)=C_c^\infty(0,\infty)$.
Next let $\sigma_+\in C_c^\infty(0,\infty)$ be the function defined in Section~\ref{S2} and 
define   the  extension $\widetilde H_+$ of $H_+$ by setting $D(\widetilde H_+)=D(H_+)+\spann \sigma_+$
and $\widetilde H_+(\varphi +\beta\,\sigma_+)=H_+\varphi-\beta\,(c\,\sigma_+')'$ for all $\varphi\in D(H_+)$ and $\beta\in\Ri$.
An analogous $L_2$-extension was used in the proof of Theorem~\ref{text4.1}.\ref{text4.1-1}.

\begin{lemma}\label{lext5.1}
The operator $\widetilde H_+$ is both $L_1$-dissipative and $L_1$-dispersive.
Therefore $ \widetilde H_+$ is $L_1$-closable and its closure  is $L_1$-dissipative and $L_1$-dispersive.
\end{lemma}
\proof\
The operator  $\widetilde H_+$  is $L_1$-dissipative if and only if 
\[
(\widetilde H_+\Phi, \sgn(\Phi))\geq 0
\]
for all $\Phi\in D(\widetilde H_+)$ where $\sgn(\Phi)$ denotes the usual sign function.
Moreover, it is $L_1$-dispersive if and only if 
\[
(\widetilde H_+\Phi, (\sgn(\Phi)\vee 0))\geq 0
\]
for all $\Phi\in D(\widetilde H_+)$.
(For background on dissipative and dispersive operators see
\cite{BaR1}, Section~2.1, or \cite{Nag}.)

Let $\Phi=\varphi+\beta\,\sigma_+$ with $\varphi\in D(H_+)$ and $\beta\in\Ri$
and note that  $(c\,\Phi')(0_+)=0$.
Next  choose a monotonically increasing $C^\infty$-function $\eta$ such that 
$\eta(x)=0$ if $|x|\leq 1$ and $\eta( x)=\pm \,1$ if $\pm \,x\geq 2$.
Then  set $\Phi_n=\eta(n\,\Phi)$.
It follows that $\Phi_n\in C_c^\infty(\Ri)$ and $\Phi_n$ converges pointwise to $\sgn(\Phi)$ as $n\to\infty$.
But integrating by parts and using  $(c\,\Phi')(0_+)=0$ one has
\[
(\widetilde H_+\Phi, \Phi_n)=\int_0^\infty dx\, c(x)\,\Phi'(x)\,\Phi_n'(x)
=\int_0^\infty dx\, c(x)\,|\Phi'(x)|^2\, n\,\eta'(n\,\Phi) \geq0
\;.
\]
Therefore in  the limit $n\to\infty $ one deduces that $(\widetilde H_+\Phi,  \sgn(\Phi))\geq 0$.
Thus $\widetilde H_+$ is $L_1$-dissipative.
The proof of dispersivity is similar.

Finally it follows by general theory that a norm densely-defined dissipative operator on a Banach space is closable and that 
its closure is dissipative. 
Moreover,  if the operator is dispersive then  the closure is also  dispersive.
(See \cite{BaR1}, Theorem~2.3.1.)
\hfill$\Box$

\bigskip

Note that as  $H_+$ is a restriction of $\widetilde H_+$ it automatically inherits the dissipativity and dispersivity properties.
Thus  $H_+$ is both dissipative and dispersive on $L_1(0,\infty)$.
Therefore its $L_1$-closure $\overline H_+^{\scriptscriptstyle 1}$
generates a strongly continuous positive contraction semigroup on $L_1(0,\infty)$ 
if and only if the range of $(I+H_+)$ is $L_1$-norm dense (see, \cite{BaR1}, Corollary~2.2.2).
But this is equivalent to the statement that  if $\psi\in L_\infty(0,\infty)$ and $(\psi,(I+H_+)\varphi)=0$
for all $\varphi\in C_c^\infty(0,\infty)$  then $\psi=0$, i.e. if $(I+H^*_+)\psi=0$ in the sense of distributions
then $\psi=0$.
If, however, $\psi-(c\,\psi')'=0$ in the distributional sense it follows that $c\,\psi'$ is  locally absolutely continuous.
Then since $c>0$ on $\langle0,\infty\rangle$ it follows that $\psi$ is a $C^1$-function locally.
Therefore $\psi=(c\,\psi')'$ in the usual sense of ordinary differential equations.

The following lemma is the key to establishing the range condition.

\begin{lemma}\label{lext5.2}
Let $ c\in W^{1,1}_{\rm loc}(0,\infty)$ be strictly positive on $\langle0,\infty\rangle$.
Assume $\int_1^\infty ds\,s\,c(s)^{-1}=\infty$.
Consider the ordinary differential equation $(c\,\psi')'=\psi$ on $\langle0,\infty\rangle$
with the boundary condition $\psi'(0_+)=\gamma\,\psi(0_+)$.

If $\gamma\geq 0$ then
there are no non-zero $L_p$-solutions $\psi$ for any $p\in[1,\infty]$.
\end{lemma}
\proof\
Assume that $\psi$ is a  non-zero solution.
Then
 \begin{eqnarray*}
 -\int^x_0|\psi|^2
 =-\int^x_0(c\,\psi')'\,\psi&=&(c\,\psi'\,\psi)(0_+)-(c\,\psi'\,\psi)(x)+\int^x_0c\,|\psi'|^2\\[5pt]
 &=&\gamma\,\psi^2(0_+)-(c\,\psi'\,\psi)(x)+\int^x_0c\,|\psi'|^2
 \;.
 \end{eqnarray*}
 Therefore
 \[
 2^{-1}(c\,(\psi^2)')(x)=\gamma\,\psi^2(0_+)+\int^x_0|\psi|^2+\int^x_0c\,|\psi'|^2\geq \int^x_0|\psi|^2\;.
 \]
 where the last bound uses $\gamma\geq 0$.
 In particular $\psi^2$ is non-decreasing.
Since  $\psi\neq0$ there is an $x_0$ such that $\psi(x_0)\neq0$ and it follows that 
  \[
 2^{-1}(c\,(\psi^2)')(x)\geq (x-x_0)\,\psi^2(x_0)
 \]
for all $x\geq x_0$.
Therefore
\[
(\psi^2)'(x)\geq \psi^2(x_0)\,x\,c(x)^{-1}
\]
for $x\geq 2\,x_0$.
It follows 
 by integration that  $|\psi(x)|^2\to\infty$ as $x\to\infty$.
Therefore there are no non-zero $L_p$-solutions.
\hfill$\Box$

\begin{remarkn}\label{rext5.1}
The conclusion of the lemma is valid in the limiting case $\gamma=+\infty$,
i.e.\ with the Dirichlet boundary condition $\psi(0_+)=0$.
\end{remarkn}

If $\psi\in L_\infty(0,\infty)$ and $(\psi,(I+H_+)\varphi)=0$ for all $\varphi\in C_c^\infty(0,\infty)$ then $\psi$
satisfies the differential equation of   Lemma~\ref{lext5.2} but it is not clear 
that it satisfies an appropriate boundary condition.
This will follow from an  $L_1$-version of Lemma~\ref{lext3.137}.

\begin{prop}\label{l3.2.3}
If $\nu_+\not\in L_\infty(0,1)$ then $\sigma_+\in D(\overline H_+^{\scriptscriptstyle 1})$
and $\overline H_+^{\scriptscriptstyle 1}\,\sigma_+=-(c\,\sigma_+')'\in L_1(0,\infty)$.
\end{prop}
 \proof\
First observe that $C_c^2(0,\infty)\subseteq D(\overline H_+^{\scriptscriptstyle 1} )$ by straightforward estimates.
Then let $\varphi_n$ be a sequence of  $C^2$-functions satisfying
$0\leq \varphi_n\leq 1$, $\varphi_n(x)=0$ if $x\in[0,n^{-1}]$ and $\varphi_n(x)=1$ if $x\geq 1$.
It follows that  $0\leq  \varphi_n\,\sigma_+\leq 1$,  $\supp\varphi_n\,\sigma_+\in[n^{-1},2]$,   $\varphi_n\,\sigma_+\in D(\overline H_+^{\scriptscriptstyle 1})$
and  $\overline H_+^{\scriptscriptstyle 1}(\varphi_n\,\sigma_+)=-(c\,(\varphi_n\,\sigma_+)')'$.

Secondly, we construct below a particular sequence of $\varphi_n$ such that 
\begin{equation}
\lim_{n\to\infty}\varphi_n(x)=1
\label{cond1}
\end{equation}
for all $x>0$ and 
\begin{equation}
\lim_{n\to\infty}\|(c\,\varphi_n')'\|_1=0  \label{cond2}
\end{equation}
 Then it follows that\[
\lim_{n\to\infty}\|\varphi_n\,\sigma_+-\sigma_+\|_1=0\;\;\;\;\;\; {\rm and}\;\;\;\;\;\;\lim_{n\to\infty}\|\overline H_+^{\scriptscriptstyle 1}(\varphi_n\, \sigma_+)+(c\,\sigma_+')'\|_1=0
\]
and the  proposition is established.

\bigskip

\noindent{\bf Construction of the sequence $\varphi_n$}$\;$
The construction is in four  steps.

\smallskip

\noindent{\bf Step 1}$\;$ 
Define $\chi_n \colon \Ri \to [0,1]$ by (\ref{apid}).
Then set $\xi_n=(1-\chi_n)^2$.
The  $\xi_n$ are positive, increasing, differentiable and 
 $\xi_n(x)\to 1$ for all $x>0$ as $n\to\infty$.
 Moreover, on  $[n^{-1},1]$ one has $\xi'_n=-2\,\chi_n'(1-\chi_n)$  and  $\xi'_n=0$  elsewhere.
But the definition of $\chi_n$ gives
\[
\xi_n'(x)=2\,c(x)^{-1}\Big(\int^x_{n^{-1}}c^{-1}\Big)\,\nu_n^{-2}
\]
for $x\in [n^{-1},1]$.
In particular $\xi_n'(n^{-1})=0$ and $\xi_n'(1)=2\,c(1)^{-1}\nu_n^{-1}$.
Thus $\xi_n$ fails to be twice-differentiable since $\xi_n'$ is discontinuous at $x=1$.
Therefore we modify the derivative by the addition of a linear function on the interval
$[n^{-1},1]$.

\smallskip

\noindent{\bf Step 2}$\;$ Define $\eta_n$ by
\begin{equation}
\eta_n(x)
= \left\{ \begin{array}{ll}
  0 & \mbox{if } x\in[0,n^{-1}\rangle \;\;\; ,  \\[5pt]
\Big(\xi_n'(x)-\xi_n'(1)(x-n^{-1})(1-n^{-1})^{-1}\Big) & \mbox{if } x\in[n^{-1}, 1] \;\;\; , \\[5pt]
  0 & \mbox{if } x \geq 1\;\;\; .\label{eex1}
         \end{array} \right.
\end{equation}
then $\eta_n(n^{-1})=0=\eta_n(1)$ and  $\eta_n$ is continuous.
Therefore setting $\zeta_n(x)=\int^x_0\eta_n$ for $x\leq 1$ and $\zeta_n(x)=\zeta_n(1)$ if $x\geq 1$
the resulting function is twice-differentiable and $\zeta_n(x)=0$ for $x\in[0,1/n]$.
Nevertheless
\[
\zeta_n(1)=\int^1_0\eta_n< \xi_n(1)=1
\]
so to complete the construction we rescale $\zeta_n$.

\smallskip

\noindent{\bf Step 3}$\;$ 
Define $\varphi_n=\zeta_n(1)^{-1}\zeta_n$.
It follows immediately that $\varphi_n(x)=0$ if $x\in[0,n^{-1}]$, $\varphi_n(x)=1$ if $x\geq 1$ and $\varphi_n$ is twice differentiable.
Moreover, 
\begin{eqnarray*}
\varphi_n(x)&=& \zeta_n(1)^{-1}\Big(1-\xi_n'(1)\int^x_{n^{-1}}ds\,(s-n^{-1})(1-n^{-1})^{-1}\Big)\\[5pt]
&\geq&\zeta_n(1)^{-1}(1-\xi_n'(1))\geq (1-\xi_n'(1))
\end{eqnarray*}
and since $\xi_n'(1)=c(1)^{-1}\nu_n^{-1}\to0$ as $n\to\infty$ one has $\varphi_n\geq 0$ for all sufficiently large $n$.
It remains to verify (\ref{cond1}) and  (\ref{cond2}).
\smallskip

\noindent{\bf Step 4}$\;$ 
First one has $1>\zeta_n(1)\geq 1-\xi_n'(1))\to 1$ as $n\to\infty$.
Therefore $\lim_{n\to\infty}\varphi_n(x)=\lim_{n\to\infty}\zeta_n(x)$.
But 
\[
\zeta_n(x)=\xi_n(x)-\xi_n'(1)\int^x_{n^{-1}}ds\,(s-n^{-1})(1-n^{-1})^{-1}\to 1
\]
for $x>0$ since $\xi_n'(1)\to0$.
Thus (\ref{cond1}) is verified.

Secondly, if $x\in[n^{-1},1]$ then 
\[
(c\,\varphi_n')(x)=\zeta_n(1)^{-1}(c\,\eta)(x)=
\zeta_n(1)^{-1}\Big((c\,\xi_n')(x)-\xi'_n(1)c(x)(x-n^{-1})(1-n^{-1})^{-1}\Big)
\;.
\]
Therefore
\[
(c\,\varphi_n')'(x)=
\zeta_n(1)^{-1}\Big((c\,\xi_n')'(x)-\xi'_n(1)c'(x)(x-n^{-1})(1-n^{-1})^{-1}-
\xi'_n(1)c(x)(1-n^{-1})^{-1}\Big)
\;.
\]
But $\zeta_n(1)^{-1}\to 1$, $\xi'_n(1)\to0$ and $\|c'\|_\infty<\infty$.
Thus
\[
\lim_{n\to\infty}\|(c\,\varphi_n')'\|_1=\lim_{n\to\infty}\|(c\,\xi_n')'\|_1
\;,
\]
i.e.\ the modifications to $\xi_n$ in Steps~2 and 3 do not affect the $L_1$-limit.
But 
\[
(c\,\xi_n')'(x)=c(x)^{-1}\,\nu_n^{-2}
\]
for $x\in[n^{-1},1]$ and is zero elsewhere.
Therefore 
\[
\|(c\,\xi_n')'\|_1=\Big(\int^1_{n^{-1}}c^{-1}\Big)\,\nu_n^{-2}=\nu_n^{-1}
\]
and  (\ref{cond2}) is verified.

This completes the proof of Proposition~\ref{l3.2.3}.\hfill$\Box$

\begin{cor}\label{cext5.1}
If $\nu_+\not\in L_\infty(0,1)$ then the $L_1$-closures of $H_+$ and $\widetilde H_+$ are equal.
\end{cor}
\proof\
This follows immediately because $D(H_+)\subseteq D(\widetilde H_+)$  but the proposition  establishes that $D(\widetilde H_+)\subseteq D(\overline H_+^{\scriptscriptstyle 1})$.
\hfill$\Box$

\begin{remarkn}\label{rext5.20}
It follows from  the proof of Proposition~\ref{l3.2.3} that if $\nu_+\not\in L_\infty(0,1)$ then one may construct a sequence $\sigma_n\in C_c^\infty(0,\infty)$ 
such that $\|\sigma_n-\sigma_+\|_1\to0$ and moreover $\|H_+\sigma_n+(c\,\sigma_+')'\|_1\to0$ as $n\to\infty$.
It suffices to replace the  $C_c^2$-approximants  $\varphi_n$ by $C_c^\infty$-approximants and to set $\sigma_n=\varphi_n\sigma_+$.
\end{remarkn}

\begin{remarkn}\label{rext5.2}
Although the foregoing results were established for $H_+=H|_{C_c^\infty(0,\infty)}$ similar statements are true for 
$H_-=H|_{C_c^\infty(-\infty,0)}$, e.g.\ $H_-$ is $L_1$-dissipative on $L_1(-\infty,0)$.
\end{remarkn}

Now we turn to the proof of Theorem~\ref{cext1.0}.
Since this involves the action of $H$ on $L_1$ and on $L_2$ it is necessary to adopt the earlier stronger assumption that $c\in W^{1,2}_{\rm loc}(\Ri)$.
Then all the preceding results apply.

\smallskip

\noindent{\bf Proof of Theorem~\ref{cext1.0}}$\;$ 
\ref{cext1.0-3}$\Rightarrow$\ref{cext1.0-2}.$\;$This is evident.

\smallskip

\noindent\ref{cext1.0-2}$\Rightarrow$\ref{cext1.0-1}.$\;$
First, it follows from an obvious extension of Lemma~\ref{lext5.1}
that $H$ is both $L_1$-dissipative and $L_1$-dispersive.
Therefore Condition~\ref{cext1.0-2} implies that the $L_1$-closure $\overline H^{\scriptscriptstyle 1}$ of $H$
generates a strongly continuous positive contraction semigroup $S$ on $L_1(\Ri)$.
Next let $H_F$ denote the Friedrichs extension of $H$ and $H_1$ the generator of the corresponding
submarkovian semigroup acting on $L_1(\Ri)$.
If $\varphi\in C_c^\infty(\Ri)$ then $\overline H^{\scriptscriptstyle 1}\varphi=H\varphi=H_F\varphi$.
But $H\varphi\in L_1(\Ri)$.
Therefore $\varphi\in D(H_1)$ and $H_1\varphi=H_F\varphi=\overline H^{\scriptscriptstyle 1}\varphi$.
Since $C_c^\infty(\Ri)$ is a core of $\overline H^1$ it follows that $H_1\supseteq \overline H^{\scriptscriptstyle 1}$.
But  $\overline H^{\scriptscriptstyle 1}$ generates a contraction semigroup and $H_1$ is $L_1$-dissipative.
Therefore $H_1=\overline H^{\scriptscriptstyle 1}$.
In particular $S$ extends to a submarkovian semigroup on the $L_p$-spaces which coincides with the 
submarkovian semigroup generated by $H_F$.

Secondly, let $\widetilde H$ be another submarkovian extension of $H$ and $\widetilde H_1$ the generator
of the corresponding semigroup on $L_1(\Ri)$.
If $\varphi\in C_c^\infty(\Ri)$  one then has $\overline H^{\scriptscriptstyle 1}\varphi=H\varphi=\widetilde H\varphi\in L_1(\Ri)$.
Therefore $\widetilde H_1\varphi=\widetilde H\varphi=\overline H^{\scriptscriptstyle 1}\varphi$ and it follows by the previous argument 
that $\widetilde H_1=\overline H^{\scriptscriptstyle 1}$.
Hence $\widetilde H_1=H_1$ and $\widetilde H$ must be the Friedrichs extension of $H$.
Therefore the Friedrichs extension is the unique submarkovian extension of $H$.

\smallskip

\ref{cext1.0-1}$\Rightarrow$\ref{cext1.0-3}.
It follows from Theorem~\ref{cext1.0} that $H$ has a unique submarkovian extension, the Friedrichs extension,
if and only if $\nu_+\vee\nu_-\not\in L_\infty(0,1)$.
Let us assume $\nu_+\not\in L_\infty(0,1)$ but $\nu_-\in L_\infty(0,1)$.
The other cases are handled similarly.

First, we argue that $(I+H_+)C_c^\infty(0,\infty)$ is dense in $L_1(0,\infty)$.
Let $\psi\in L_\infty(0,\infty)$ such that $(\psi,(I+H_+)\varphi)=0$ for all $\varphi\in C_c^\infty(0,\infty)$. 
Then $\psi$ satisfies the ordinary differential equation $(c\,\psi')'=\psi$.
 In particular $ c\,\psi'$ is continuous near the origin and $(c\,\psi')(0_+)$ exists.
 But if $(c\,\psi')(0_+)\neq0$ then $\psi\sim \nu_+$ as $x\to0_+$.
 Since $\nu_+\not\in L_\infty(0,1)$ this contradicts the boundedness of $\psi$.
 Therefore $(c\,\psi')(0_+)=0$.
Then the assumed growth conditions at infinity allow the application of  Lemma~\ref{lext5.2} and one  deduces that $\psi=0$.
 Therefore $(I+H_+)C_c^\infty(0,\infty)$ is dense in $L_1(0,\infty)$.

Secondly, consider the restriction $H_-$ of $H$ to $C_c^\infty(-\infty,0)$.
Since $\nu_-\in L_\infty(0,1)$ one cannot apply the foregoing reasoning
to establish that $(I+H_-)C_c^\infty(-\infty,0)$ is dense in $L_1(-\infty,0)$.
It follows, however, from Proposition~\ref{l3.2.3} that $\sigma_+\in D(\overline H_+^{\scriptscriptstyle 1})$, or $\sigma_+\in D(\overline H^{\scriptscriptstyle 1})$ on $L_1(\Ri)$.
Since $\sigma_+\in C_c^\infty(\Ri)=D(H)$ it follows that $\sigma_-=\sigma-\sigma_+\in D(\overline H^{\scriptscriptstyle 1})$.
Now one can define a dissipative extension $\widetilde H_-$ of $H_-$ by
 $D(\widetilde H_-)=D(H_-)+\spann \sigma_-$
and $\widetilde H_-(\varphi +\beta\,\sigma_-)=H_-\varphi-\beta\,(c\,\sigma_-')'$ for all $\varphi\in D(H_-)$ and $\beta\in\Ri$.
Then $\overline H\supseteq \widetilde H_-\oplus H_+$ and to deduce that $(I+H)C_c^\infty(\Ri\backslash\{0\})$ is dense in $L_1(\Ri)$
it  suffices to prove that $(I+\widetilde H_-)D(\widetilde H_-)$ is dense in $L_1(-\infty,0)$.
 
Fourthly, suppose there is a $\psi\in L_\infty(-\infty,0)$ such that $(\psi, (I+\widetilde H_-)\Phi)=0$ for all $\Phi=\varphi +\beta\,\sigma_-\in D(\widetilde H_-)$.
It follows that 
\[
(\psi, (I+ H_-)\varphi)+\beta\,\Big((\psi,\sigma_-)-(\psi, (c\,\sigma_-')')\Big)=0
\]
for all $\varphi\in C_c^\infty(-\infty,0)$ and all $\beta\in \Ri$.
Therefore $(\psi, (I+ H_-)\varphi)=0$ for all $\varphi\in C_c^\infty(-\infty,0)$ as before.
In addition, however, one must have 
\[
(\psi,\mu)-(\psi,(c\,\mu')')=0
\]
for all $\mu=\sigma_- \;(\mod C_c^\infty(-\infty,0))$.
But integration by parts gives 
\[
(\psi,\mu)=((c\,\psi')',\mu)-(\mu\,c\,\psi')(0_-)
\]
again for all $\mu=\sigma_- \;(\mod C_c^\infty(-\infty,0))$.
This immediately implies that  $(c\,\psi')(0_-)=0$.
Finally, arguing as above,  there is no non-zero bounded $\psi$ satisfying $(\psi, (I+ H_-)\varphi)=0$
for all $\varphi\in C_c^\infty(-\infty,0)$ and $(c\,\psi')(0_-)=0$.
So $(I+\widetilde H_-)D(\widetilde H_-)$ is dense in $L_1(-\infty,0)$.\hfill$\Box$

\bigskip

Conditions \ref{cext1.0-3} and \ref{cext1.0-2} of Theorem~\ref{cext1.0} imply that $C_c^\infty(\Ri)$ is a core
of the generator $H_1$ of the submarkovian semigroup $S$ acting on $L_1(\Ri)$.
This in turn implies that the semigroup is conservative. 
Indeed one has $(H\varphi,\one)=0$ for all $\varphi\in C_c^\infty(\Ri)$ and  by closure for all $\varphi\in D(H_1)$.
Then, however,
\[
{{d}\over{dt}}(S_t\varphi,\one)=(H_1S_t\varphi,\one)=0
\]
$\varphi\in D(H_1)$.
Therefore $S$ must be conservative.
Davies \cite{Dav14}, Theorem~2.2, has established a converse statement for a large class of elliptic operators on $\Ri^d$
(see also \cite{Pan2}).

\section{Lipschitz coefficients}\label{S4}

In this section we examine operators with  a coefficient  $c\in W^{1,\infty}_{\rm loc}(\Ri)$
and give a proof of Theorem~\ref{cext1.1}.
The simplest case is for $c$ strictly positive on $\Ri\backslash\{0\}$ but $c(0)=0$.
Then the $W^{1,\infty}_{\rm loc}$-assumption on $c$ ensures that $c(x)=O(x)$ as $x\to0_\pm$.
Thus $\nu_\pm\not\in L_\infty(0,1)$ and $H$ has a unique submarkovian extension by 
Theorem~\ref{text4.1}.\ref{text4.1-2}.
Another simple situation occurs if  $c(x)=0$ for $x\leq 0$ but $c(x)>0$ for $x>0$. 
Then the Lipschitz condition
means that $c(x)=O(x)$ as $x\to0_+$ and the uniqueness follows from 
Theorem~\ref{tnsm1.1}.\ref{tnsm1.1-2}.
To understand the general situation one needs information about the extensions of $H$
acting on a finite interval with the coefficient degenerate at both endpoints.
These extensions have been extensively studied by Feller  \cite{Feller} \cite{Feller2} \cite{Feller3}
(see also \cite{Mandl}) for $H$  acting on the spaces $C_b$ and $L_1$ using probabilistic arguments and by 
  Ulmet  using function analytic techniques \cite{Ulm}.
Properties of the Friedrichs extension on $L_2$  have also been analyzed  in detail by Campiti, Metafune and 
Pallara  \cite{CMP}.  
But all  self-adjoint extensions can also be studied by the   methods of the previous sections.
The situation for the submarkovian extensions is particularly simple.

Define  $\nu(x)=\int^{1/2}_xc^{-1}$.
Fix   $\sigma_0\in C_c^\infty(0,1/2)$ with $0\leq \sigma_0\leq 1$ and  $\sigma_0=1$ in a neighbourhood of zero.
Then set $\tau_0=\sigma_0\,\nu$.
Further define $\sigma_1$ and $\tau_1$ as the reflections of $\sigma_0$ and $\tau_0$ around the midpoint $1/2$ of the interval.

\begin{lemma}\label{lext5.10}
Assume $c\in W^{1,\infty}(0,1)$, $c>0$ on $\langle0,1\rangle$ and $c(0)=0=c(1)$.
Define the symmetric operator $H$ on $L_2(0,1)$ by $H\varphi=-(c\,\varphi')'$   for  $\varphi\in C_c^\infty(0,1)$.

Then  $H$ has a unique submarkovian extension  $H_F$, the Friedrichs extension, 
and $D(H_F)=D(\overline H)+\spann\sigma_0+\spann\sigma_1$.
\end{lemma}
\proof\
Since $c\in W^{1,\infty}(0,1)$ one has $c(x)=O(x)$ as $x\to0_+$ and $c(x)=O(1-x)$ as $x\to1_-$.
Therefore $\nu$ is unbounded at both  endpoints $0$ and $1$.
Then the proof follows the arguments used for the half-line in Section~\ref{S2}.
In fact $ H_F(\varphi+a_0\,\sigma_0+a_1\sigma_1)=\overline H\varphi -a_0(c\,\sigma_0')'
-a_1(c\,\sigma_1')'$ for all $\varphi\in D(\overline H)$ and $a_0,a_1\in\Ri$
in  direct analogy with  Remark~\ref{rext2.13} for the half-line and Remark~\ref{rext3.11}
for the line.
\hfill$\Box$

\bigskip

It also follows under the assumptions of the lemma that $(I+H)C_c^\infty(0,1)$ is dense in $L_1(0,1)$.
This can be deduced from the argument used in the proof of Lemma~\ref{lext5.2} or from Proposition~3.5 of \cite{CMP}.

\smallskip

The foregoing lemma is the last  element in the proof of  Theorem~\ref{cext1.1}.

\smallskip

\noindent{\bf Proof of Theorem~\ref{cext1.1}}$\;$
First, let $\cz=\{x\in\Ri: c(x)=0\}$ denote the zero set of the coefficient $c$.
Then $L_2(\Ri)=L_2(\cz)\oplus L_2(\cz^{\rm c})$.
Since $H=0$ in restriction to
$L_2(\cz)$ we only need to analyze the operator on $ L_2(\cz^{\rm c})$.
Next  as $c$ is continuous  $\cz$ is closed and the complement $\cz^{\rm c}$ is  open.
Therefore $\cz^{\rm c}$ is the union of a family of  disjoint open intervals $I_i$.

Secondly, $H|_{C_c^\infty(I_i)}$ has a unique submarkovian extension $H_i$ by Theorem~\ref{tnsm1.1}  if $I_i$ is semi-infinite and by Lemma~\ref{lext5.10} 
if $I_i$ is finite.
In both cases $H_i$ is the Friedrichs extension and its action is given either by the algorithm of Remark~\ref{rext2.13}
or by that of Lemma~\ref{lext5.10}. 
In particular if $I_i=\langle a_i,\infty\rangle$ then $D(H_i)=D(\overline H^{\scriptscriptstyle i})+\spann\sigma_{a_i}$ where $\overline H^{\scriptscriptstyle i}$ denotes the $L_2(I_i)$-closure of 
$H|_{C_c^\infty(I_i)}$ and $\sigma_{a_i}$ is a $C_c^\infty(I_i)$-function which is equal to one in a neighbourhood of~$a_i$.
Similarly if $I_i=\langle a_i,b_i\rangle$ then $D(H_i)=D(\overline H^{\scriptscriptstyle i})+\spann\sigma_{a_i}+\spann \sigma_{b_i}$ 
and if   $I_i=\langle -\infty,b_i\rangle$ then $D(H_i)=D(\overline H^{\scriptscriptstyle i})+\spann \sigma_{b_i}$.

Thirdly, set $\widetilde H=\bigoplus_iH_i$.
Then 
$\widetilde H$ is a submarkovian extension of $H|_{C_c^\infty(\cz^{\rm c})}$
corresponding to the Friedrichs extension of $H|_{C_c^\infty(\cz^{\rm c})}$.
In particular $D(\widetilde H)=\bigoplus_iD(H_i)$ consists of the $\varphi=\bigoplus_i\varphi_i$ with $\varphi_i\in D(H_i)$ 
such that $\sum_i(\|\varphi_i\|_{L_2(I_i)}^2+\|H_i\varphi_i\|_{L_2(I_i)}^2)<\infty$ and then 
$\widetilde H\varphi=\bigoplus_iH_i\varphi_i$.
Now let $\widehat H$ denote a second submarkovian extension of $H$.
We will prove that $\widehat H|_{D(H_i)}=H_i$ for each $i$ and thereby deduce that $\widehat H\supseteq \widetilde H$.
But a self-adjoint operator cannot have a proper self-adjoint extension so one must have $\widehat H=\widetilde H$,
i.e.\  $H$ has a unique submarkovian extension.

If  $\varphi\in C_c^\infty(I_i)\subseteq D(H)$ then $\widehat H\varphi=H\varphi=\widetilde H\varphi$.
Therefore  $\widehat H\varphi=\widetilde H\varphi$ for all $\varphi\in D(\overline H^{\scriptscriptstyle i})$.
Next suppose $I_i$ has a finite left endpoint $a_i$ and $\sigma_{a_i}\in D(H_i)$ but $\sigma_{a_i}\not\in D(\overline H^{\scriptscriptstyle i})$.
It follows from the proof of Proposition~\ref{l3.2.3} (see Remark~\ref{rext5.20}) that one may choose a sequence $\sigma_{n,i}\in C_c^\infty(I_i)$
such that $\|\sigma_{n,i} -\sigma_{a_i}\|_1\to0$ and $\|H\sigma_{n,i}+(c\,\sigma_{a_i}')'\|_1\to0$ as $n\to\infty$.
Now let  $\widehat H_1$ and $\widetilde H_1$ denote the $L_1$-generators of the submarkovian semigroups generated by $\widehat H$ and $\widetilde H$, respectively.
Then $\widehat H_1\sigma_{n,i}=H\sigma_{n,i}=\widetilde H_1\sigma_{n,i}$.
But $\widehat H_1$ is $L_1$-closed. 
Therefore $\sigma_{a_i}\in D(\widehat H_1)$ and $\widehat H_1\sigma_{a_i}=-(c\,\sigma_{a_i}')'$.
In addition~$\sigma_{a_i}$ and $(c\,\sigma_{a_i}')'$ are both in $L_2(I_i)$.
Hence $\sigma_{a_i}\in D(\widehat H_1)\cap L_2(I_i)$ and $\widehat H_1\sigma_{a_i}=-(c\,\sigma_{a_i}')'\in L_2(I_i)$.
Therefore  $\sigma_{a_i}\in D(\widehat H)$ and  $\widehat H\sigma_{a_i}= \widehat H_1\sigma_{a_i}=-(c\,\sigma_{a_i}')'$.
Similarly $\sigma_{a_i}\in D(\widetilde  H)$ and $\widetilde H\sigma_{a_i}=-(c\,\sigma_{a_i}')'$.
Thus one concludes that  $\widehat H\sigma_{a_i}=\widetilde H\sigma_{a_i}$.
If $I_i$ has a finite right endpoint~$b_i$ one concludes similarly that $\widehat H\sigma_{b_i}=\widetilde H\sigma_{b_i}$.
Therefore $\widehat H$ and $\widetilde H$ are equal on $D(H_i)$.
This completes the proof of uniqueness of the submarkovian extension.\hfill$\Box$

\subsection*{Acknowledgement}
This project originated in a series of interesting discussions with Ricardo Weder and Gian Michele Graf whilst the first author was visiting 
the Institut f\"ur Theoretische Physik at the Eidgen\"ossische Technische Hochschule, Z\"urich in 2007.
The author is indebted to J\"urg Fr\"ohlich and Gian Michele Graf  for facilitating this visit and to the ETH for providing financial support.

\end{document}